\documentclass[draft]{elsart}

\usepackage{amsmath}
\usepackage{amsfonts}
\usepackage{amssymb}
\usepackage{graphicx}
\usepackage[all]{xy}

\newtheorem{introcor}{Corollary}

\newtheorem{introthm}{Theorem}

\numberwithin{equation}{section}
\numberwithin{thm}{section}

\newcommand{\BBr}{\mathop{\mathcal{B}\mathrm{r}}\nolimits}

\begin{document}

\begin{frontmatter}



\title{Mackey-functor structure on the Brauer groups of a finite Galois covering of schemes}


\author{Hiroyuki NAKAOKA\thanksref{JSPS}\thanksref{Prof}}

\address{Graduate School of Mathematical Sciences, The University of Tokyo 
3-8-1 Komaba, Meguro, Tokyo, 153-8914 Japan}

\ead{deutsche@ms.u-tokyo.ac.jp}
\thanks[JSPS]{The author is supported by JSPS.}
\thanks[Prof]{The author wishes to thank Professor Toshiyuki Katsura for his encouragement.}

\begin{abstract}

Past studies of the Brauer group of a scheme tells us the importance of the interrelationship among Brauer groups of its finite \'{e}tale coverings.
In this paper, we consider these groups simultaneously, and construct an integrated object \lq\lq Brauer-Mackey functor".

We realize this as a {\it cohomological Mackey functor} on the Galois category of finite \'{e}tale coverings.
For any finite \'{e}tale covering of schemes, we can associate two homomorphisms for Brauer groups, namely the pull-back and the norm map. These homomorphisms make Brauer groups into a bivariant functor ($=$ Mackey functor) on the Galois category.

As a corollary, Restricting to a finite Galois covering of schemes, we obtain a cohomological Mackey functor on its Galois group. This is a generalization of the result for rings by Ford \cite{Ford}. Moreover, applying Bley and Boltje's theorem \cite{Bley-Boltje}, we can derive certain isomorphisms for the Brauer groups of intermediate coverings.
\end{abstract}

\end{frontmatter}


\section{Introduction}
\label{Introduction}
In this paper, any scheme $X$ is assumed to be Noetherian. 
$\pi(X)$ denotes its \'{e}tale fundamental group. 
Any morphism is locally of finite type, unless otherwise specified. 
As in \cite{Milne}, $X_{\mathrm{et}}$ denotes the small \'{e}tale site, consisting of \'{e}tale morphisms of finite type over $X$. 
If $\mathcal{U}=(U_i\overset{f_i}{\longrightarrow}X)_{i\in I}$ is a covering in this site, we write as $\mathcal{U}\in\mathrm{Cov}_{\mathrm{et}}(X)$. 
$\mathcal{U}\prec\mathcal{V}$ means $\mathcal{U}$ is a refinement of $\mathcal{V}$.

 As for the finite \'{e}tale covering, the \'{e}tale fundamental group and the Galois category, we follow the terminology in \cite{Murre}. For example a finite \'{e}tale covering is just a finite \'{e}tale morphism of schemes.

Our aim is to make the following generalization of the result by Ford \cite{Ford}, which was shown for rings.
\begin{introcor}[Corollary \ref{MainCor1}]
Let $\pi:Y\rightarrow X$ be a finite Galois covering of schemes with Galois group $G$. Assume $X$ satisfies Assumption \ref{Assumption}. Then the correspondence
\[ H\le G\mapsto \mathrm{Br}(Y/H) \]
forms a cohomological Mackey functor on $G$. Here, $H\le G$ means $H$ is a subgroup of $G$.
\end{introcor}

This follows from our main theorem:
\begin{introthm}[Theorem \ref{MainThm}]
Let $S$ be a connected scheme satisfying Assumption \ref{Assumption}. Let $(\mathrm{FEt}/S)$ denote the category of finite \'{e}tale coverings over $S$. Then, the Brauer group functor $\mathrm{Br}$ forms a cohomological Mackey functor on $(\mathrm{FEt}/S)$.
\end{introthm}

As in Definition \ref{DefOfMackFtr}, a Mackey functor is a bivariant pair of functors $\mathrm{Br}=(\mathrm{Br}^{\ast},\mathrm{Br}_{\ast})$. For any morphism $\pi:Y\rightarrow X$, the contravariant part $\mathrm{Br}^{\ast}(\pi):\mathrm{Br}(X)\rightarrow \mathrm{Br}(Y)$ is the pull-back, and the covariant part $\mathrm{Br}_{\ast}(\pi):\mathrm{Br}(Y)\rightarrow \mathrm{Br}(X)$ is the norm map defined later.

By applying Bley and Boltje's theorem (Fact \ref{Bley-Boltje's Thm}) to Corollary \ref{MainCor1}, we can obtain certain relations between Brauer groups of intermediate coverings:
\begin{introcor}[Corollary \ref{MainCor2}]
Let $X$ be a connected scheme satisfying Assumption \ref{Assumption}, and $\pi:Y\rightarrow X$ be a finite Galois covering with $\mathrm{Gal}(Y/X)=G$.\\
{\rm (i)} Let $\ell$ be a prime number. If $H\le G$ is not $\ell$-hypoelementary, then there is a natural isomorphism of $\mathbb{Z}_{\ell}$-modules
\[ \bigoplus _{\underset{n:\text{odd}}{U=H_0<\cdots<H_n=H}}\mathrm{Br}(Y/U)(\ell)^{|U|}\cong\bigoplus_{\underset{n:\text{even}}{U=H_0<\cdots<H_n=H}}\mathrm{Br}(Y/U)(\ell)^{|U|}. \]\\
{\rm (ii)} If $H\le G$ is not hypoelementary, then there is a natural isomorphism of abelian groups
\[ \bigoplus _{\underset{n:\text{odd}}{U=H_0<\cdots<H_n=H}}\mathrm{Br}(Y/U)^{|U|}\cong\bigoplus_{\underset{n:\text{even}}{U=H_0<\cdots<H_n=H}}\mathrm{Br}(Y/U)^{|U|}. \]
Here, $|U|$ denotes the order of $U$.
\end{introcor}

\section{Preliminaries}
\label{Preliminaries}

To fix the notation, we recall several facts in this section. 
If $\mathcal{C}$ is a category and $X$ is an object in $\mathcal{C}$, we abbreviately write as $X\in\mathcal{C}$. 
If $f:X\rightarrow Y$ is a morphism in $\mathcal{C}$, we write as $f\in\mathcal{C}(X,Y)$ or $f\in\mathrm{Mor}_{\mathcal{C}}(X,Y)$. 

Monoidal categories, monoidal functors and monoidal transformations are always assumed to be symmetric.

For a scheme $X$, $\mathrm{q\text{-}Coh}(X)$ denotes the category of quasi-coherent modules over $\mathcal{O}_X$.

\subsection{Fpqc descent}

\begin{defn}\label{DescentDef}
Let $X^{\prime}\rightarrow X$ be an fpqc morphism of schemes.
Put $X^{(2)}:=X^{\prime}\times_XX^{\prime}$, $X^{(3)}:=X^{\prime}\times_XX^{\prime}\times_XX^{\prime}$ and let
\begin{align*}
p_i&:X^{(2)}\rightarrow X^{\prime}\qquad (i=1,2)\\
p_{ij}&:X^{(3)}\rightarrow X^{(2)}\quad\ (i,j\in\{1,2,3\})
\end{align*}
be the projections. Define a category $\mathrm{q\text{-}Coh}(X^{\prime}\rightarrow X)$ as follows $:$

- an object in $\mathrm{q\text{-}Coh}(X^{\prime}\rightarrow X)$ is a pair $(\mathcal{F},\varphi)$ of a sheaf $\mathcal{F}\in\mathrm{q\text{-}Coh}(X^{\prime})$ and an isomorphism $\varphi:p_1^{\ast}\mathcal{F}\overset{\cong}{\longrightarrow}p_2^{\ast}\mathcal{F}$ in $\mathrm{q\text{-}Coh}(X^{(2)})$.

- a morphism from $(\mathcal{F},\varphi)$ to $(\mathcal{G},\psi)$ is a morphism $\alpha\in\mathrm{q\text{-}Coh}(X^{\prime})(\mathcal{F},\mathcal{G})$, such that
\[ p_2^{\ast}\alpha\circ\phi=\psi\circ p_1^{\ast}\alpha. \]

For any $(\mathcal{F},\varphi)$ and $(\mathcal{G},\psi)\in\mathrm{q\text{-}Coh}(X^{\prime}\rightarrow X)$, let $\varphi\otimes\psi$ be the abbreviation of
\begin{align*}
p_1^{\ast}(\mathcal{F}\underset{\mathcal{O}_{X^{\prime}}}{\otimes}\mathcal{G})%
\overset{\cong}{\rightarrow}%
p_1^{\ast}\mathcal{F}\underset{\mathcal{O}_{X^{(2)}}}{\otimes}p_1^{\ast}\mathcal{G}%
\overset{\varphi\otimes\psi}{\longrightarrow}%
p_2^{\ast}\mathcal{F}\underset{\mathcal{O}_{X^{(2)}}}{\otimes}p_2^{\ast}\mathcal{G}%
\overset{\cong}{\rightarrow}%
p_2^{\ast}(\mathcal{F}\underset{\mathcal{O}_{X^{\prime}}}{\otimes}\mathcal{G}).
\end{align*}
Then, $\mathrm{q\text{-}Coh}(X^{\prime}\rightarrow X)$ has a canonical symmetric monoidal structure defined by
\[ (\mathcal{F},\varphi)\otimes(\mathcal{G},\psi):=(\mathcal{F}\underset{\mathcal{O}_{X^{\prime}}}{\otimes}\mathcal{G},\varphi\otimes\psi). \]
\end{defn}

\begin{rem}\label{DescentRem}
Let $f:X^{\prime}\rightarrow X$ be an fpqc morphism of schemes.
The pull-back functor by $f$
\[ f^{\ast}:\mathrm{q\text{-}Coh}(X)\rightarrow\mathrm{q\text{-}Coh}(X^{\prime}) \]
factors through $\mathrm{q\text{-}Coh}(X^{\prime}\rightarrow X)$ $:$
\[
\xy
(0,8)*+{\mathrm{q\text{-}Coh}(X^{\prime}\rightarrow X)}="0";
(0,-14)*+{}="1";
(-32,-6)*+{\mathrm{q\text{-}Coh}(X)}="2";
(32,-6)*+{\mathrm{q\text{-}Coh}(X^{\prime})}="4";
{\ar^{\exists\underline{f^{\ast}}} "2";"0"};
{\ar^{U} "0";"4"};
{\ar@/_0.80pc/_{f^{\ast}} "2";"4"};
{\ar@{}|\circlearrowright "0";"1"};
\endxy
\]
where $U$ is the forgetful functor.
By the fpqc descent, $\underline{f^{\ast}}$ is an equivalence.

In fact, $U$ is a monoidal functor, and $\underline{f^{\ast}}$ is a monoidal equivalence.
\end{rem}

\subsection{Contravariant nature of the Brauer group}

\begin{rem}
Let $\pi:Y\rightarrow X$ be a finite \'{e}tale covering. 
For any abelian sheaf $\mathcal{G}$ on $Y_{\mathrm{et}}$ and any positive integer $q$, the following composition of the canonical morphisms is an isomorphism $:$
\[ \mathfrak{c}:H^q_{\mathrm{et}}(X,\pi_{\ast}\mathcal{G})\rightarrow H^q_{\mathrm{et}}(Y,\pi^{\ast}\pi_{\ast}\mathcal{G})\rightarrow H^q_{\mathrm{et}}(Y,\mathcal{G}) \]
\end{rem}

\begin{rem}
For any scheme $X$, there exists a natural monomorphism
\[ \chi_X:\mathrm{Br}(X)\hookrightarrow \mathrm{Br}^{\prime}(X):=H^2_{\mathrm{et}}(X,\mathbb{G}_{m,X})_{\mathrm{tor}}\ , \]
such that for any morphism $\pi:Y\rightarrow X$, 
\[
\xy
(-16,7)*+{\mathrm{Br}(X)}="0";
(-16,-7)*+{H^2_{\mathrm{et}}(X,\mathbb{G}_{m,X})}="2";
(16,7)*+{\mathrm{Br}(Y)}="4";
(16,-7)*+{H^2_{\mathrm{et}}(Y,\mathbb{G}_{m,Y})}="6";
(-16,3)*+{}="10";
(16,3)*+{}="14";
{\ar@{^(->}^<<<<<<{\chi_X} "10";"2"};
{\ar@{^(->}^<<<<<<{\chi_Y} "14";"6"};
{\ar^{\pi^{\ast}} "0";"4"};
{\ar^{\pi^{\ast}} "2";"6"};
{\ar@{}|\circlearrowright "0";"6"};
\endxy
\]
is a commutative diagram.
\end{rem}
Here $\pi^{\ast}:\mathrm{Br}(X)\rightarrow \mathrm{Br}(Y)$ is the pull-back of Azumaya algebras, while $\pi^{\ast}:H^2_{\mathrm{et}}(X,\mathbb{G}_{m,X})\rightarrow H^2_{\mathrm{et}}(Y,\mathbb{G}_{m,Y})$ is defined as the composition of the canonical morphism $\mathfrak{c}:H^2_{\mathrm{et}}(X,\pi_{\ast}\mathbb{G}_{m,Y})\rightarrow H^2_{\mathrm{et}}(Y,\mathbb{G}_{m,Y})$ and $H^2_{\mathrm{et}}(\pi_{\sharp}):H^2_{\mathrm{et}}(X,\mathbb{G}_{m,X})\rightarrow H^2_{\mathrm{et}}(X,\pi_{\ast}\mathbb{G}_{m,Y})$, where $\pi_{\sharp}:\mathbb{G}_{m,X}\rightarrow \pi_{\ast}\mathbb{G}_{m,Y}$ is the canonical (structure) homomorphism of \'{e}tale sheaves on $X$.

\section{Norm functor}

In this section, we construct a monoidal functor
\[ \mathcal{N}_{\pi}:\mathrm{q\text{-}Coh}(Y)\rightarrow \mathrm{q\text{-}Coh}(X) \]
which we call the {\it norm functor}, for any finite \'{e}tale covering $\pi:Y\rightarrow X$.

\subsection{Trivial case}
\begin{defn}\label{TrivNormDef}
Let $X$ be a scheme, and let
\[ \nabla=\nabla_{X,d}:\underset{1\le k\le d}{\textstyle{\coprod}}X_k\rightarrow X\quad (X_k=X\ (1\le\forall k\le d)) \]
be the folding map. We define the norm functor
\[ \mathcal{N}_{\nabla}:\mathrm{q\text{-}Coh}(\underset{1\le k\le d}{\textstyle{\coprod}}X_k)\rightarrow \mathrm{q\text{-}Coh}(X) \]
by
\[ \mathcal{N}_{\nabla}(\mathcal{G}):=\mathcal{G}\mid_{X_1}\otimes_{\mathcal{O}_X}\cdots \otimes_{\mathcal{O}_X}\mathcal{G}\mid_{X_d} \]
for any $\mathcal{G}\in \mathrm{q\text{-}Coh}(\underset{1\le k\le d}{\textstyle{\coprod}}X_k)$, and similarly for morphisms.
\end{defn}

\begin{rem}
$\mathcal{N}_{\nabla}$ is a monoidal functor.
\end{rem}

\begin{rem}\label{DisjRem}
For any automorphism $\tau:\underset{1\le k\le d}{\coprod}X_k\overset{\cong}{\longrightarrow}\underset{1\le k\le d}{\coprod}X_k$ compatible with $\nabla$, there is a natural monoidal isomorphism
\[ \mathcal{N}_{\nabla}\circ\tau^{\ast}\cong\mathcal{N}_{\nabla}. \]
\end{rem}
\begin{pf}
Left to the reader.
\end{pf}

\begin{defn}\label{TrivDef}
Let $\pi:Y\rightarrow X$ be a finite \'{e}tale covering. Assume there exists an isomorphism
\[ \eta:\underset{1\le k\le d}{\textstyle{\coprod}}X_k\overset{\cong}{\longrightarrow}Y \]
compatible with $\pi$ and $\nabla_{X,d}$. We define the norm functor $\mathcal{N}_{\pi}$ by
\[ \mathcal{N}_{\pi}:=\mathcal{N}_{\nabla}\circ\eta^{\ast}.\]
\end{defn}

\begin{rem}
By Remark \ref{DisjRem}, $\mathcal{N}_{\pi}$ does not depend on the choice of trivialization $\eta$.
\end{rem}

\begin{rem}\label{TrivializationRem}
Let $\pi:Y\rightarrow X$ be a finite \'{e}tale covering with a trivialization $\eta:\underset{1\le k\le d}{\textstyle{\coprod}}X_k\overset{\cong}{\longrightarrow}Y$, as in Definition \ref{TrivDef}. Let $f:X^{\prime}\rightarrow X$ be any morphism and take the pull-back $:$
\begin{align}
\xy
(-10,8)*+{Y^{\prime}}="0";
(10,8)*+{Y}="2";
(-10,-8)*+{X^{\prime}}="4";
(10,-8)*+{X}="6";
(0,0)*+{\square}="8";
{\ar^{g} "0";"2"};
{\ar_{\pi^{\prime}} "0";"4"};
{\ar_{f} "4";"6"};
{\ar^{\pi} "2";"6"};
\endxy
\label{PullBackDiag1}
\end{align}

Then by pulling $\eta$ back by $f$, we obtain an isomorphism
\[ \eta^{\prime}:\underset{1\le k\le d}{\textstyle{\coprod}}X^{\prime}_k\overset{\cong}{\longrightarrow}Y^{\prime},\quad (X^{\prime}_k=X^{\prime}\ (1\le\forall k\le d)) \]
compatible with $\pi^{\prime}$ and $\nabla_{X^{\prime},d}$ $:$
\[
\xy
(-14,12)*+{Y^{\prime}}="0";
(-14,-1.6)*+{}="1";
(14,12)*+{Y}="2";
(-14,-2.4)*+{}="3";
(-14,-12)*+{X^{\prime}}="4";
(14,-12)*+{X}="6";
(-24,-2)*+{\underset{d}{\textstyle{\coprod}}X^{\prime}}="8";
(4,-2)*+{\underset{d}{\textstyle{\coprod}}X}="10";
{\ar^{g} "0";"2"};
{\ar@{-} "0";"1"};
{\ar^{\pi^{\prime}} "3";"4"};
{\ar_{f} "4";"6"};
{\ar^{\pi} "2";"6"};
{\ar^{\eta^{\prime}} "8";"0"};
{\ar^{\eta} "10";"2"};
{\ar_{\nabla} "8";"4"};
{\ar_{\nabla} "10";"6"};
{\ar^>>>>>>>>>{\underset{d}{\textstyle{\coprod}}f} "8";"10"};
\endxy
\]
\begin{center}
$($ all faces are commutative $)$
\end{center}
\end{rem}

\begin{prop}\label{TrivthetaProp}
Let $\pi:Y\rightarrow X$ be a finite \'{e}tale covering with a trivialization.

{\rm (i)} For any morphism $f:X^{\prime}\rightarrow X$, if we take the pull-back as in $(\ref{PullBackDiag1})$, then there exists a natural monoidal isomorphism
\[ \theta^f_{\pi}:f^{\ast}\circ\mathcal{N}_{\pi}\overset{\cong}{\longrightarrow}\mathcal{N}_{\pi^{\prime}}\circ g^{\ast}. \]

Moreover, $\theta$ is natural in $f$ $:$

{\rm (ii)} For any other morphism $f^{\prime}:X^{\prime\prime}\rightarrow X^{\prime}$, if we take the pull-back
\[
\xy
(-10,8)*+{Y^{\prime\prime}}="0";
(10,8)*+{Y^{\prime}}="2";
(-10,-8)*+{X^{\prime\prime}}="4";
(10,-8)*+{X^{\prime}}="6";
(0,0)*+{\square}="8";
(13,-9)*+{,}="10";
{\ar^{g^{\prime}} "0";"2"};
{\ar_{\pi^{\prime\prime}} "0";"4"};
{\ar_{f^{\prime}} "4";"6"};
{\ar^{\pi^{\prime}} "2";"6"};
\endxy
\]
then we have
\begin{equation}
\theta^{f\circ f^{\prime}}_{\pi}=(\theta^{f^{\prime}}_{\pi^{\prime}}\circ g^{\ast})\cdot(f^{\prime\ast}\circ\theta^f_{\pi}).
\label{Equationtheta}
\end{equation}
\[
\xy
(-20,10)*+{(f\circ f^{\prime})^{\ast}\mathcal{N}_{\pi}}="1";
(20,10)*+{\mathcal{N}_{\pi^{\prime\prime}}(g\circ g^{\prime})^{\ast}}="2";
(-28,-2)*+{f^{\prime\ast}f^{\ast}\mathcal{N}_{\pi}}="3";
(28,-2)*+{\mathcal{N}_{\pi^{\prime\prime}}g^{\prime\ast}g^{\ast}}="4";
(0,-12)*+{f^{\prime\ast}\mathcal{N}_{\pi^{\prime}}g^{\ast}}="5";
(0,12)*+{}="6";
{\ar^{\theta^{f\circ f^{\prime}}_{\pi}} "1";"2"};
{\ar^{\cong} "1";"3"};
{\ar_{\cong} "2";"4"};
{\ar_{f^{\prime\ast}\circ\theta^f_{\pi}} "3";"5"};
{\ar_{\theta^{f^{\prime}}_{\pi^{\prime}}\circ g^{\ast}} "5";"4"};
{\ar@{}|\circlearrowright "5";"6"};
\endxy
\]
\end{prop}
\begin{pf}
{\rm (i)} This follows from Remark \ref{TrivializationRem}, since we have
\begin{eqnarray*}
f^{\ast}\mathcal{N}_{\pi}(\mathcal{G})&=& f^{\ast}((\eta^{\ast}\mathcal{G})\mid_{X_1}\otimes_{\mathcal{O}_X}\cdots\otimes_{\mathcal{O}_X}(\eta^{\ast}\mathcal{G})\mid_{X_d})\\
&\cong&f^{\ast}((\eta^{\ast}\mathcal{G})\mid_{X_1})\otimes_{\mathcal{O}_{X^{\prime}}}\cdots\otimes_{\mathcal{O}_{X^{\prime}}}f^{\ast}((\eta^{\ast}\mathcal{G})\mid_{X_d})\\
&\cong&((\underset{d}{\textstyle{\coprod}}f)^{\ast}\eta^{\ast}\mathcal{G})\mid_{X^{\prime}_1}\otimes_{\mathcal{O}_{X^{\prime}}}\cdots\otimes_{\mathcal{O}_{X^{\prime}}}((\underset{d}{\textstyle{\coprod}}f)^{\ast}\eta^{\ast}\mathcal{G})\mid_{X^{\prime}_d}\\
&\cong&(\eta^{\prime\ast}(g^{\ast}\mathcal{G}))\mid_{X^{\prime}_1}\otimes_{\mathcal{O}_{X^{\prime}}}\cdots\otimes_{\mathcal{O}_{X^{\prime}}}(\eta^{\prime\ast}(g^{\ast}\mathcal{G}))\mid_{X^{\prime}_d}\\
&=&\mathcal{N}_{\pi^{\prime}}(g^{\ast}\mathcal{G})
\end{eqnarray*}
for any $\mathcal{G}\in\mathrm{q\text{-}Coh}(Y)$, in the notation of Remark \ref{TrivializationRem}

{\rm (ii)}
This follows from the trivial case
\begin{align*}
&\xy
(-32,8)*+{\underset{1\le k\le d}{\textstyle{\coprod}}X^{\prime\prime}_k}="0";
(0,8)*+{\underset{1\le k\le d}{\textstyle{\coprod}}X^{\prime}_k}="2";
(32,8)*+{\underset{1\le k\le d}{\textstyle{\coprod}}X_k}="4";
(-32,-8)*+{X^{\prime\prime}}="10";
(0,-8)*+{X^{\prime}}="12";
(32,-8)*+{X}="14";
(-16,0)*+{\square}="1";
(16,0)*+{\square}="3";
{\ar^{g^{\prime}} "0";"2"};
{\ar^{g} "2";"4"};
{\ar_{f^{\prime}} "10";"12"};
{\ar_{f} "12";"14"};
{\ar_{\nabla_{X^{\prime\prime},d}} "0";"10"};
{\ar_{\nabla_{X^{\prime},d}} "2";"12"};
{\ar^{\nabla_{X,d}} "4";"14"};
\endxy\\
&\qquad\qquad\qquad\qquad(g=\underset{d}{\textstyle{\coprod}}f,\ g^{\prime}=\underset{d}{\textstyle{\coprod}}f^{\prime}),
\end{align*}
where, $(\ref{Equationtheta})$ immediately follows from the commutativity of the following diagrams for any $\mathcal{F}\in\mathrm{q\text{-}Coh}(Y)\,:$
\[
\xy
(-32,10)*+{(f\circ f^{\prime})^{\ast}(\underset{1\le k\le d}{\bigotimes}\mathcal{F}\mid_{X_k})}="1";
(32,10)*+{\underset{1\le k\le d}{\bigotimes}((f\circ f^{\prime})^{\ast}(\mathcal{F}\mid_{X_k}))}="2";
(-40,-6)*+{f^{\prime\ast}(f^{\ast}(\underset{1\le k\le d}{\bigotimes}\mathcal{F}\mid_{X_k}))}="3";
(40,-6)*+{\underset{1\le k\le d}{\bigotimes}(f^{\prime\ast}(f^{\ast}(\mathcal{F}\mid_{X_k})))}="4";
(0,-20)*+{f^{\prime\ast}(\underset{1\le k\le d}{\bigotimes}(f^{\ast}(\mathcal{F}\mid_{X_k})))}="5";
(0,16)*+{}="6";
{\ar^{\cong} "1";"2"};
{\ar^{\cong} "1";"3"};
{\ar_{\cong} "2";"4"};
{\ar_{\cong} "3";"5"};
{\ar_{\cong} "5";"4"};
{\ar@{}|\circlearrowright "5";"6"};
\endxy
\]
\[
\xy
(-20,10)*+{(f\circ f^{\prime})^{\ast}(\mathcal{F}\mid_{X_k})}="1";
(20,10)*+{((g\circ g^{\prime})^{\ast}\mathcal{F})\mid_{X^{\prime\prime}_k}}="2";
(-28,-2)*+{f^{\prime\ast}(f^{\ast}(\mathcal{F}\mid_{X_k}))}="3";
(28,-2)*+{(g^{\prime\ast}(g^{\ast}\mathcal{F}))\mid_{X^{\prime\prime}_k}}="4";
(0,-12)*+{f^{\prime\ast}((g^{\ast}\mathcal{F})\mid_{X^{\prime}_k})}="5";
(0,12)*+{}="6";
{\ar^{\cong} "1";"2"};
{\ar^{\cong} "1";"3"};
{\ar_{\cong} "2";"4"};
{\ar_{\cong} "3";"5"};
{\ar_{\cong} "5";"4"};
{\ar@{}|\circlearrowright "5";"6"};
\endxy
\quad(1\le \forall k\le d)
\]
\end{pf}

\subsection{Constant degree case}

\begin{rem}\label{LocalTrivRem}
Let $\pi:Y\rightarrow X$ be a finite \'{e}tale covering of constant degree $d$. There exists an fpqc morphism $f:X^{\prime}\rightarrow X$ such that the pull-back of $\pi$ by $f$ becomes trivial $:$
\[
\xy
(-10,8)*+{Y^{\prime}}="0";
(-24,0)*+{\underset{d}{\textstyle{\coprod}}X^{\prime}}="1";
(10,8)*+{Y}="2";
(-10,-8)*+{X^{\prime}}="4";
(10,-8)*+{X}="6";
(0,0)*+{\square}="8";
(-4,0)*+{}="10";
{\ar^{g} "0";"2"};
{\ar^{\cong} "1";"0"};
{\ar_{\nabla_{X^{\prime},d}} "1";"4"};
{\ar^{\pi^{\prime}} "0";"4"};
{\ar_{f} "4";"6"};
{\ar^{\pi} "2";"6"};
{\ar@{}|\circlearrowright "1";"10"};
\endxy
\]

$f$ can be also taken as a surjective \'{e}tale morphism.
\end{rem}

\begin{prop}
In the notation of Remark \ref{LocalTrivRem},
\[ \mathcal{N}_{\pi^{\prime}}\circ g^{\ast}:\mathrm{q\text{-}Coh}(Y)\rightarrow \mathrm{q\text{-}Coh}(X^{\prime}) \]
factors through $\mathrm{q\text{-}Coh}(X^{\prime}\overset{f}{\rightarrow}X)\overset{U}{\longrightarrow}\mathrm{q\text{-}Coh}(X^{\prime})$ in Remark \ref{DescentRem}.
\end{prop}
\begin{pf}
For the convenience, we abbreviate two functors
\begin{align*}
g^{\ast}&:\mathrm{q\text{-}Coh}(Y)\rightarrow\mathrm{q\text{-}Coh}(Y^{\prime})\\
\mathcal{N}_{\pi^{\prime}}\circ g^{\ast}&:\mathrm{q\text{-}Coh}(Y)\rightarrow\mathrm{q\text{-}Coh}(X^{\prime})
\end{align*}
respectively to
\begin{align*}
\widetilde{\ \cdot\ }&:\mathrm{q\text{-}Coh}(Y)\rightarrow\mathrm{q\text{-}Coh}(Y^{\prime})\\
\overline{\ \cdot\ }&:\mathrm{q\text{-}Coh}(Y)\rightarrow\mathrm{q\text{-}Coh}(X^{\prime}).
\end{align*}
These are monoidal. Put
\begin{align*}
X^{(2)}:=X^{\prime}\times_XX^{\prime},&\ \ Y^{(2)}:=Y^{\prime}\times_YY^{\prime},\\
X^{(3)}:=X^{\prime}\times_XX^{\prime}\times_XX^{\prime},&\ \ Y^{(3)}:=Y^{\prime}\times_YY^{\prime}\times_YY^{\prime},
\end{align*}
and denote the projections by
\begin{align*}
p_i:X^{(2)}\rightarrow X^{\prime},&\ \ q_i:Y^{(2)}\rightarrow Y^{\prime}\quad(i=1,2),\\
p_{ij}:X^{(3)}\rightarrow X^{(2)},&\ \ q_{ij}:Y^{(3)}\rightarrow Y^{(2)} \quad(1\le i<j\le 3),\\
p^{\prime}_i:X^{(3)}\rightarrow X^{\prime},&\ \ q^{\prime}_i:Y^{(3)}\rightarrow Y^{\prime}\quad(i=1,2,3).
\end{align*}

Pulling $\pi$ back by these projections, we obtain finite \'{e}tale coverings:
\begin{align*}
\pi^{(2)}:Y^{(2)}\rightarrow X^{(2)}\\
\pi^{(3)}:Y^{(3)}\rightarrow X^{(3)}
\end{align*}
\[
\xy
(-30,8)*+{Y^{(3)}}="10";
(-10,8)*+{Y^{(2)}}="12";
(-30,-8)*+{X^{(3)}}="20";
(-10,-8)*+{X^{(2)}}="22";
(10,8)*+{Y^{\prime}}="14";
(30,8)*+{Y}="16";
(10,-8)*+{X^{\prime}}="24";
(30,-8)*+{X}="26";
(-20,0)*+{\square}="1";
(0,0)*+{\square}="2";
(20,0)*+{\square}="3";
{\ar^{q_{ij}} "10";"12"};
{\ar_{\pi^{(3)}} "10";"20"};
{\ar_{p_{ij}} "20";"22"};
{\ar^{\pi^{(2)}} "12";"22"};
{\ar^{q_{\ell}} "12";"14"};
{\ar_{p_{\ell}} "22";"24"};
{\ar^{g} "14";"16"};
{\ar_{\pi^{\prime}} "14";"24"};
{\ar_{f} "24";"26"};
{\ar^{\pi} "16";"26"};
{\ar^{g} "14";"16"};
{\ar_{\pi^{\prime}} "14";"24"};
{\ar_{f} "24";"26"};
{\ar^{\pi} "16";"26"};
\endxy
\]
\begin{center}
$(1\le i<j\le 3, \ell=1,2)$
\end{center}

Remark that each of $\pi^{(2)}$ and $\pi^{(3)}$ has a trivialization.

It suffices to show the following:
\begin{claim}\label{Claim11}
For each $\mathcal{F}\in \mathrm{q\text{-}Coh}(Y)$, there is a canonical isomorphism
\[ \phi_{\mathcal{F}}:p_1^{\ast}\overline{\mathcal{F}}\overset{\cong}{\longrightarrow} p_2^{\ast}\overline{\mathcal{F}} \]
such that for any morphism $\alpha\in\mathrm{q\text{-}Coh}(Y)(\mathcal{F},\mathcal{G})$,
\[ p_2^{\ast}\alpha\circ\phi_{\mathcal{F}}=\phi_{\mathcal{G}}\circ p_1^{\ast}\alpha \]
is satisfied.
\end{claim}
\begin{pf}(proof of Claim \ref{Claim11})
Since $\widetilde{\mathcal{F}}=g^{\ast}\mathcal{F}$, there is a canonical isomorphism $\psi_{\mathcal{F}}:q_1^{\ast}\widetilde{\mathcal{F}}\overset{\cong}{\longrightarrow}q_2^{\ast}\widetilde{\mathcal{F}}$
such that
\begin{align}
q_{13}^{\ast}\psi_{\mathcal{F}}=q_{23}^{\ast}\psi_{\mathcal{F}}\circ q_{12}^{\ast}\psi_{\mathcal{F}}.
\label{Equationpsi}
\end{align}

Define $\phi_{\mathcal{F}}:p_1^{\ast}\overline{\mathcal{F}}\overset{\cong}{\rightarrow} p_2^{\ast}\overline{\mathcal{F}}$ by
\begin{align*}
 \phi_{\mathcal{F}}:&=(\theta^{p_2}_{\pi^{\prime}})^{-1}\circ\mathcal{N}_{\pi^{(2)}}(\psi_{\mathcal{F}})\circ\theta^{p_1}_{\pi^{\prime}}\\
&=(p_1^{\ast}\mathcal{N}_{\pi^{\prime}}\widetilde{\mathcal{F}}%
\overset{\theta^{p_1}_{\pi^{\prime}}}{\longrightarrow}%
\mathcal{N}_{\pi^{(2)}}q_1^{\ast}\widetilde{\mathcal{F}}%
\overset{\mathcal{N}_{\pi^{(2)}}(\psi_{\mathcal{F}})}{\longrightarrow}%
\mathcal{N}_{\pi^{(2)}}q_2^{\ast}\widetilde{\mathcal{F}}%
\overset{(\theta^{p_2}_{\pi^{\prime}})^{-1}}{\longrightarrow}%
p_2^{\ast}\mathcal{N}_{\pi^{\prime}}\widetilde{\mathcal{F}})
\end{align*}

By $(\ref{Equationtheta})$ and the naturality of $\theta$, we have a commutative diagram
\[
\xy
(-18,10)*+{p_i^{\prime\ast}\overline{\mathcal{F}}}="0";
(18,10)*+{p_j^{\prime\ast}\overline{\mathcal{F}}}="2";
(-18,-10)*+{\mathcal{N}_{\pi^{(3)}}q_i^{\prime\ast}\widetilde{\mathcal{F}}}="4";
(18,-10)*+{\mathcal{N}_{\pi^{(3)}}q_j^{\prime\ast}\widetilde{\mathcal{F}}}="6";
{\ar^{p_{ij}^{\ast}\phi_{\mathcal{F}}} "0";"2"};
{\ar_{\theta^{p_i^{\prime}}_{\pi^{\prime}}} "0";"4"};
{\ar_{\mathcal{N}_{\pi^{(3)}}(q_{ij}^{\ast}\psi_{\mathcal{F}})} "4";"6"};
{\ar^{\theta^{p_j^{\prime}}_{\pi^{\prime}}} "2";"6"};
{\ar@{}|\circlearrowright "0";"6"};
\endxy
\]
for each $1\le i<j\le 3$.
Thus $p_{13}^{\ast}\phi_{\mathcal{F}}=p_{23}^{\ast}\phi_{\mathcal{F}}\circ p_{12}^{\ast}\phi_{\mathcal{F}}$ follows from $(\ref{Equationpsi})$.
\end{pf}
\end{pf}

\begin{rem}\label{NormpiefRem}
For any $\mathcal{F},\mathcal{G}\in\mathrm{q\text{-}Coh}(Y)$, we have a commutative diagram
\[
\xy
(-20,14)*+{p_1^{\ast}(\overline{\mathcal{F}\otimes\mathcal{G}})}="1";
(20,14)*+{p_2^{\ast}(\overline{\mathcal{F}\otimes\mathcal{G}})}="2";
(-28,0)*+{p_1^{\ast}(\overline{\mathcal{F}}\otimes\overline{\mathcal{G}})}="3";
(28,0)*+{p_2^{\ast}(\overline{\mathcal{F}}\otimes\overline{\mathcal{G}})}="4";
(-20,-14)*+{p_1^{\ast}\overline{\mathcal{F}}\otimes p_1^{\ast}\overline{\mathcal{G}}}="5";
(20,-14)*+{p_2^{\ast}\overline{\mathcal{F}}\otimes p_2^{\ast}\overline{\mathcal{G}}}="6";
(32,-15)*+{.}="7";
{\ar^{\phi_{\mathcal{F}\otimes\mathcal{G}}} "1";"2"};
{\ar^{\cong} "1";"3"};
{\ar_{\cong} "2";"4"};
{\ar^{\cong} "3";"5"};
{\ar_{\cong} "4";"6"};
{\ar_{\phi_{\mathcal{F}}\otimes\phi_{\mathcal{G}}} "5";"6"};
{\ar@{}|\circlearrowright "1";"6"};
\endxy
\]
From this, we can see easily that the factorization of $\mathcal{N}_{\pi^{\prime}}\circ g^{\ast}$
\[ \underline{\mathcal{N}_{\pi^{\prime}}\circ g^{\ast}}:\mathrm{q\text{-}Coh}(Y)\rightarrow\mathrm{q\text{-}Coh}(X^{\prime}\rightarrow X) \]
becomes a monoidal functor.
\end{rem}

\begin{cor}
Let $\pi:Y\rightarrow X$ be a finite \'{e}tale covering of constant degree $d$ and $f:X^{\prime}\rightarrow X$ be an fpqc morphism trivializing $\pi$. Then we have a monoidal functor
\[ \mathcal{N}_{\pi}^f:\mathrm{q\text{-}Coh}(Y)\rightarrow\mathrm{q\text{-}Coh}(X)\]
uniquely up to a natural monoidal isomorphism, such that $($in the notation of Remark \ref{LocalTrivRem}$\, )$, there is a natural monoidal isomorphism
\begin{equation}\label{Moriast}
\underline{f^{\ast}}\circ\mathcal{N}^f_{\pi}\cong\underline{\mathcal{N}_{\pi^{\prime}}\circ g^{\ast}}\, ,
\end{equation}
and thus $f^{\ast}\circ\mathcal{N}^f_{\pi}\cong\mathcal{N}_{\pi^{\prime}}\circ g^{\ast}$.
\end{cor}
\begin{pf}
This follows from Remark \ref{DescentRem} and Remark \ref{NormpiefRem}.
\end{pf}

\begin{prop}\label{ConstNormProp}
Let $\pi:Y\rightarrow X$ be a finite \'{e}tale covering of constant degree $d$.
If $f_1:X_1\rightarrow X$ and $f_2:X_2\rightarrow X$ are fpqc morphisms trivializing $\pi$, then there exists a natural monoidal isomorphism
\[ \mathcal{N}^{f_1}_{\pi}\cong\mathcal{N}^{f_2}_{\pi}. \]
\end{prop}
\begin{pf}
By considering the pull-back
\[
\xy
(-20,0)*+{X_1\times_XX_2}="0";
(0,8)*+{X_1}="2";
(0,-8)*+{X_2}="4";
(16,0)*+{X}="6";
(19,-2)*+{,}="7";
(0,0)*+{\square}="8";
{\ar^{} "0";"2"};
{\ar_{} "0";"4"};
{\ar_{f_2} "4";"6"};
{\ar^{f_1} "2";"6"};
\endxy
\]
we may assume $f_2$ factors through $f_1$:
\[
\xy
(0,4)*+{X_1}="0";
(-20,-4)*+{X_2}="2";
(20,-4)*+{X}="4";
(0,-10)*+{}="6";
{\ar^{\exists f_3} "2";"0"};
{\ar^{f_1} "0";"4"};
{\ar@/_0.80pc/_{f_2} "2";"4"};
{\ar@{}|\circlearrowright "0";"6"};
\endxy
\]

Pulling $\pi$ back by $f_i$ for each $i=1,2$,  we obtain diagrams
\[
\xy
(-10,8)*+{Y_i}="0";
(-32,4)*+{\underset{1\le k\le d}{\textstyle{\coprod}}X_{i,k}}="1";
(10,8)*+{Y}="2";
(-10,-8)*+{X_i}="4";
(10,-8)*+{X}="6";
(14,-9)*+{,}="7";
(0,0)*+{\square}="8";
(-2,-2)*+{}="10";
{\ar^{g_i} "0";"2"};
{\ar^>>>>>>>>>{\eta_i}_>>>>>>>>>{\cong} "1";"0"};
{\ar_{\nabla_{X_i,d}} "1";"4"};
{\ar^{\pi_i} "0";"4"};
{\ar_{f_i} "4";"6"};
{\ar^{\pi} "2";"6"};
{\ar@{}|\circlearrowright "1";"10"};
\endxy
\qquad
\xy
(-10,8)*+{Y_2}="0";
(10,8)*+{Y_1}="2";
(-10,-8)*+{X_2}="4";
(10,-8)*+{X_1}="6";
(0,0)*+{\square}="8";
{\ar^{g_3} "0";"2"};
{\ar_{\pi_2} "0";"4"};
{\ar_{f_3} "4";"6"};
{\ar^{\pi_1} "2";"6"};
\endxy
\]
where $X_{i,k}=X_i\ (1\le \forall k\le d)$.

Let
\[
\xy
(-8,8)*+{X_1^{(2)}}="0";
(8,8)*+{X_1}="2";
(-8,-8)*+{X_1}="4";
(8,-8)*+{X}="6";
(0,0)*+{\square}="8";
{\ar^{p_2} "0";"2"};
{\ar_{p_1} "0";"4"};
{\ar_{f_1} "4";"6"};
{\ar^{f_1} "2";"6"};
\endxy
\qquad
\xy
(-8,8)*+{X_2^{(2)}}="0";
(8,8)*+{X_2}="2";
(-8,-8)*+{X_2}="4";
(8,-8)*+{X}="6";
(0,0)*+{\square}="8";
{\ar^{p_2^{\prime}} "0";"2"};
{\ar_{p_1^{\prime}} "0";"4"};
{\ar_{f_2} "4";"6"};
{\ar^{f_2} "2";"6"};
\endxy
\]
be pull-backs, and let $f_3^{(2)}:X_2^{(2)}\rightarrow X_1^{(2)}$ be the induced morphism:
\[
\xy
(-8,8)*+{X_2^{(2)}}="0";
(8,8)*+{X_1^{(2)}}="2";
(-8,-8)*+{X_2}="4";
(8,-8)*+{X_1}="6";
{\ar^{f_3^{(2)}} "0";"2"};
{\ar_{p_i^{\prime}} "0";"4"};
{\ar_{f_3} "4";"6"};
{\ar^{p_i} "2";"6"};
{\ar@{}|\circlearrowright "0";"6"};
\endxy
\qquad
(i=1,2)
\]

Using $(\ref{Equationtheta})$, we can show that the natural monoidal isomorphism
\[ \Xi_{\mathcal{F}}:=(f_3^{\ast}\mathcal{N}_{\pi_1}g_1^{\ast}\mathcal{F}\overset{\theta^{f_3}_{\pi_1}\circ g_1^{\ast}}{\longrightarrow}\mathcal{N}_{\pi_2}g_3^{\ast}g_1^{\ast}\mathcal{F}\overset{\cong}{\rightarrow}\mathcal{N}_{\pi_2}g_2^{\ast}\mathcal{F})\quad(\mathcal{F}\in\mathrm{q\text{-}Coh}(Y)) \]
is compatible with descent data
\begin{align*}
\phi_{\mathcal{F}}&: p_1^{\ast}\mathcal{N}_{\pi_1}g_1^{\ast}\mathcal{F}\overset{\cong}{\longrightarrow}p_2^{\ast}\mathcal{N}_{\pi_1}g_1^{\ast}\mathcal{F}\\
\phi^{\prime}_{\mathcal{F}}&: p_1^{\prime\ast}\mathcal{N}_{\pi_2}g_2^{\ast}\mathcal{F}\overset{\cong}{\longrightarrow}p_2^{\prime\ast}\mathcal{N}_{\pi_2}g_2^{\ast}\mathcal{F}
\end{align*}
defined in Claim \ref{Claim11}:
\[
\xy
(4,0)*+{\circlearrowright}="0";
(-8,16)*+{p_1^{\prime\ast}f_3^{\ast}\mathcal{N}_{\pi_1}g_1^{\ast}\mathcal{F}}="1";(-40,8)*+{f_3^{(2)\ast}p_1^{\ast}\mathcal{N}_{\pi_1}g_1^{\ast}\mathcal{F}}="2";
(-40,-8)*+{f_3^{(2)\ast}p_2^{\ast}\mathcal{N}_{\pi_1}g_1^{\ast}\mathcal{F}}="3";
(-8,-16)*+{p_2^{\prime\ast}f_3^{\ast}\mathcal{N}_{\pi_1}g_1^{\ast}\mathcal{F}}="4";
(30,-16)*+{p_2^{\prime\ast}\mathcal{N}_{\pi_2}g_2^{\ast}\mathcal{F}}="5";
(30,16)*+{p_1^{\prime\ast}\mathcal{N}_{\pi_2}g_2^{\ast}\mathcal{F}}="6";
{\ar^{\cong} "1";"2"};
{\ar_{f_3^{(2)\ast}\phi_{\mathcal{F}}} "2";"3"};
{\ar_{\cong} "4";"3"};
{\ar^{p^{\prime\ast}_1\Xi_{\mathcal{F}}} "1";"6"};
{\ar_{p^{\prime\ast}_2\Xi_{\mathcal{F}}} "4";"5"};
{\ar^{\phi^{\prime}_{\mathcal{F}}} "6";"5"};
\endxy
\qquad
(i=1,2)
\]
Thus Proposition \ref{ConstNormProp} follows from Remark \ref{DescentRem}.
\end{pf}

\begin{defn}
Let $\pi:Y\rightarrow X$ be a finite \'{e}tale covering of constant degree $d$. By Proposition \ref{ConstNormProp}, $\mathcal{N}^f_{\pi}$ is uniquely determined as a monoidal functor up to a natural monoidal isomorphism, independently of the fpqc morphism $f$ trivializing $\pi$. We denote this functor simply by $\mathcal{N}_{\pi}$, and call it the {\it norm functor} for $\pi$.
\end{defn}

\begin{prop}\label{ConstthetaProp}
Let $\pi:Y\rightarrow X$ be a finite \'{e}tale covering of constant degree $d$.
Let $f:X^{\prime}\rightarrow X$ be a morphism, and take the pull-back
\[
\xy
(-10,8)*+{Y^{\prime}}="0";
(10,8)*+{Y}="2";
(-10,-8)*+{X^{\prime}}="4";
(10,-8)*+{X}="6";
(0,0)*+{\square}="8";
(14,-9)*+{.}="10";
{\ar^{g} "0";"2"};
{\ar_{\pi^{\prime}} "0";"4"};
{\ar_{f} "4";"6"};
{\ar^{\pi} "2";"6"};
\endxy
\]
Then there exists a natural monoidal isomorphism
\[ \theta^f_{\pi}:f^{\ast}\mathcal{N}_{\pi}\overset{\cong}{\longrightarrow}\mathcal{N}_{\pi^{\prime}}g^{\ast}. \]
\end{prop}
\begin{pf}
Let $u:U\rightarrow X$ be an fpqc morphism trivializing $\pi$, and take the pull-backs
\[
\xy
(-8,8)*+{V}="0";
(8,8)*+{Y}="2";
(-8,-8)*+{U}="4";
(8,-8)*+{X}="6";
(0,0)*+{\square}="8";
(11,-9)*+{,}="10";
{\ar^{v} "0";"2"};
{\ar_{\varpi} "0";"4"};
{\ar_{u} "4";"6"};
{\ar^{\pi} "2";"6"};
\endxy
\quad
\xy
(-8,8)*+{U^{\prime}}="0";
(8,8)*+{U}="2";
(-8,-8)*+{X^{\prime}}="4";
(8,-8)*+{X}="6";
(0,0)*+{\square}="8";
(11,-9)*+{,}="10";
{\ar^{f_U} "0";"2"};
{\ar_{u^{\prime}} "0";"4"};
{\ar_{f} "4";"6"};
{\ar^{u} "2";"6"};
\endxy
\quad
\xy
(-8,8)*+{V^{\prime}}="0";
(8,8)*+{V}="2";
(-8,-8)*+{Y^{\prime}}="4";
(8,-8)*+{Y}="6";
(0,0)*+{\square}="8";
(11,-9)*+{,}="10";
{\ar^{g_V} "0";"2"};
{\ar_{v^{\prime}} "0";"4"};
{\ar_{g} "4";"6"};
{\ar^{v} "2";"6"};
\endxy
\quad
\xy
(-8,8)*+{V^{\prime}}="0";
(8,8)*+{Y^{\prime}}="2";
(-8,-8)*+{U^{\prime}}="4";
(8,-8)*+{X^{\prime}}="6";
(0,0)*+{\square}="8";
(11,-9)*+{.}="10";
{\ar^{v^{\prime}} "0";"2"};
{\ar_{\varpi^{\prime}} "0";"4"};
{\ar_{u^{\prime}} "4";"6"};
{\ar^{\pi^{\prime}} "2";"6"};
\endxy
\]
Remark
\[
\xy
(-8,8)*+{V^{\prime}}="0";
(8,8)*+{V}="2";
(-8,-8)*+{U^{\prime}}="4";
(8,-8)*+{U}="6";
(0,0)*+{\square}="8";
{\ar^{g_V} "0";"2"};
{\ar_{\varpi^{\prime}} "0";"4"};
{\ar_{f_U} "4";"6"};
{\ar^{\varpi} "2";"6"};
\endxy
\]
is also a pull-back diagram.

By Proposition \ref{TrivthetaProp}, there is a natural monoidal isomorphism
\[ \theta^{f_U}_{\varpi}:f_U^{\ast}\mathcal{N}_{\varpi}\overset{\cong}{\longrightarrow}\mathcal{N}_{\varpi^{\prime}}g_V^{\ast}. \]
As in Proposition \ref{ConstNormProp}, natural monoidal isomorphism
\begin{equation*}
\Theta^f_{\pi}:f_U^{\ast}\mathcal{N}_{\varpi}v^{\ast}
\overset{\theta^{f_U}_{\varpi}\circ v^{\ast}}{\longrightarrow}%
\mathcal{N}_{\varpi^{\prime}}g_V^{\ast}v^{\ast}%
\overset{\cong}{\rightarrow}%
\mathcal{N}_{\varpi^{\prime}}v^{\prime\ast}g^{\ast}%
\end{equation*}
is compatible with descent data, and we obtain a natural monoidal isomorphism
\[ \theta^f_{\pi}:f^{\ast}\mathcal{N}_{\pi}\overset{\cong}{\longrightarrow}\mathcal{N}_{\pi^{\prime}}g^{\ast} \]
such that $u^{\prime\ast}\theta^f_{\pi}$ gives $\Theta^f_{\pi}$.
\end{pf}

As in Proposition \ref{TrivthetaProp}, $\theta$ is natural in $f\,:$
\begin{cor}\label{ConstthetaCor}
Let $\pi:Y\rightarrow X$ be a finite \'{e}tale covering of constant degree.
For any morphisms $X^{\prime\prime}\overset{f^{\prime}}{\longrightarrow}X^{\prime}\overset{f}{\longrightarrow}X$, if we take the pull-back
\[
\xy
(-16,8)*+{Y^{\prime\prime}}="10";
(0,8)*+{Y^{\prime}}="12";
(16,8)*+{Y}="14";
(-16,-8)*+{X^{\prime\prime}}="20";
(0,-8)*+{X^{\prime}}="22";
(16,-8)*+{X}="24";
(-8,0)*+{\square}="2";
(8,0)*+{\square}="4";
(19,-9)*+{,}="6";
{\ar^{g^{\prime}} "10";"12"};
{\ar_{\pi^{\prime\prime}} "10";"20"};
{\ar_{f^{\prime}} "20";"22"};
{\ar^{\pi^{\prime}} "12";"22"};
{\ar^{g} "12";"14"};
{\ar_{f} "22";"24"};
{\ar^{\pi} "14";"24"};
\endxy
\]
then we have
\[ \theta^{f\circ f^{\prime}}_{\pi}=(\theta^{f^{\prime}}_{\pi^{\prime}}\circ g^{\ast})\cdot(f^{\prime\ast}\circ\theta^f_{\pi}). \]
\[
\xy
(-20,10)*+{(f\circ f^{\prime})^{\ast}\mathcal{N}_{\pi}}="1";
(20,10)*+{\mathcal{N}_{\pi^{\prime\prime}}(g\circ g^{\prime})^{\ast}}="2";
(-28,-2)*+{f^{\prime\ast}f^{\ast}\mathcal{N}_{\pi}}="3";
(28,-2)*+{\mathcal{N}_{\pi^{\prime\prime}}g^{\prime\ast}g^{\ast}}="4";
(0,-12)*+{f^{\prime\ast}\mathcal{N}_{\pi^{\prime}}g^{\ast}}="5";
(0,12)*+{}="6";
{\ar^{\theta^{f\circ f^{\prime}}_{\pi}} "1";"2"};
{\ar^{\cong} "1";"3"};
{\ar_{\cong} "2";"4"};
{\ar_{f^{\prime\ast}\circ\theta^f_{\pi}} "3";"5"};
{\ar_{\theta^{f^{\prime}}_{\pi^{\prime}}\circ g^{\ast}} "5";"4"};
{\ar@{}|\circlearrowright "5";"6"};
\endxy
\]
\end{cor}
\begin{pf}
Let $u:U\rightarrow X$ be an fpqc morphism trivializing $\pi$, and take the pull-backs
\[
\xy
(-8,8)*+{V}="0";
(8,8)*+{Y}="2";
(-8,-8)*+{U}="4";
(8,-8)*+{X}="6";
(0,0)*+{\square}="8";
{\ar^{v} "0";"2"};
{\ar_{\varpi} "0";"4"};
{\ar_{u} "4";"6"};
{\ar^{\pi} "2";"6"};
\endxy
\qquad
\xy
(-16,12)*+{V^{\prime\prime}}="0";
(0,12)*+{V^{\prime}}="2";
(16,12)*+{V}="4";
(-16,0)*+{U^{\prime\prime}}="10";
(0,0)*+{U^{\prime}}="12";
(16,0)*+{U}="14";
(-16,-12)*+{X^{\prime\prime}}="20";
(0,-12)*+{X^{\prime}}="22";
(16,-12)*+{X}="24";
(-8,6)*+{\square}="1";
(8,6)*+{\square}="3";
(-8,-6)*+{\square}="11";
(8,-6)*+{\square}="13";
(19,-13)*+{,}="6";
{\ar_{\varpi^{\prime\prime}} "0";"10"};
{\ar^{\varpi^{\prime}} "2";"12"};
{\ar^{\varpi} "4";"14"};
{\ar_{u^{\prime\prime}} "10";"20"};
{\ar^{u^{\prime}} "12";"22"};
{\ar^{u} "14";"24"};
{\ar^{g^{\prime}_{V^{\prime}}} "0";"2"};
{\ar^{g_V} "2";"4"};
{\ar_{f^{\prime}_{U^{\prime}}} "10";"12"};
{\ar_{f_U} "12";"14"};
{\ar_{f^{\prime}} "20";"22"};
{\ar_{f} "22";"24"};
\endxy
\]
Applying Proposition \ref{TrivthetaProp}, we obtain the following commutative diagram:
\[
\xy
(-20,10)*+{(f_U\circ f_{U^{\prime}}^{\prime})^{\ast}\mathcal{N}_{\varpi}v^{\ast}}="1";
(20,10)*+{\mathcal{N}_{\varpi^{\prime\prime}}v^{\prime\prime\ast}(g\circ g^{\prime})^{\ast}}="2";
(-28,-2)*+{f_{U^{\prime}}^{\prime\ast}f_U^{\ast}\mathcal{N}_{\varpi}v^{\ast}}="3";
(28,-2)*+{\mathcal{N}_{\varpi^{\prime\prime}}v^{\prime\prime\ast}g^{\prime\ast}g^{\ast}}="4";
(0,-16)*+{f^{\prime\ast}_{U^{\prime}}\mathcal{N}_{\varpi^{\prime}}v^{\prime\ast}g^{\ast}}="5";
(0,12)*+{}="6";
{\ar^{\Theta^{f\circ f^{\prime}}_{\pi}} "1";"2"};
{\ar^{\cong} "1";"3"};
{\ar_{\cong} "2";"4"};
{\ar_{f^{\prime\ast}_{U^{\prime}}\circ\Theta^f_{\pi}} "3";"5"};
{\ar_{\Theta^{f^{\prime}}_{\pi^{\prime}}\circ g^{\ast}} "5";"4"};
{\ar@{}|\circlearrowright "5";"6"};
\endxy
\]
From this, we obtain
\[ u^{\prime\prime\ast}\circ\theta^{f\circ f^{\prime}}_{\pi}=(u^{\prime\prime\ast}\circ\theta^{f^{\prime}}_{\pi^{\prime}}\circ g^{\ast})\cdot(u^{\prime\prime\ast}\circ f^{\prime\ast}\circ\theta^f_{\pi}). \]
Since $u^{\prime\prime}$ is fpqc, Corollary \ref{ConstthetaCor} follows.
\end{pf}

\subsection{General case}

\begin{rem}
Let $X$ be a scheme. For any open subscheme $\iota:U\hookrightarrow X$ and $\mathcal{H}\in\mathrm{q\text{-}Coh}(U)$, we often abbreviate $\iota_{\ast}\mathcal{H}\in\mathrm{q\text{-}Coh}(X)$ simply to $\mathcal{H}$.

Let $X=\underset{1\le i\le n}{\coprod}X_i$ be the decomposition into the connected open components. For any $\mathcal{F}\in\mathrm{q\text{-}Coh}(X)$, we have a canonical decomposition
\[ \mathcal{F}=\underset{1\le i\le n}{\bigoplus}\mathcal{F}\mid_{X_i}=\mathcal{F}\mid_{X_1}\oplus\cdots\oplus\mathcal{F}\mid_{X_n}. \]
Regarding this decomposition, for any $\mathcal{F}, \mathcal{G}\in\mathrm{q\text{-}Coh}(X)$, we have
\[ \mathcal{F}\underset{\mathcal{O}_X}{\otimes}\mathcal{G}=(\mathcal{F}\mid_{X_1}\underset{\mathcal{O}_X}{\otimes}\mathcal{G}\mid_{X_1})\oplus\cdots\oplus(\mathcal{F}\mid_{X_n}\underset{\mathcal{O}_X}{\otimes}\mathcal{G}\mid_{X_n}). \]
\end{rem}

\begin{defn}
Let $\pi:Y\rightarrow X$ be a finite \'{e}tale covering, and let $X=\underset{1\le i\le n}{\coprod}X_i$ be the decomposition into the connected open components.
Put $Y_i:=\pi^{-1}(X_i)$, and let $\pi_i:Y_i\rightarrow X_i$ be the restriction of $\pi$ onto $Y_i$.

We define the {\it norm functor}
\[ \mathcal{N}_{\pi}:\mathrm{q\text{-}Coh}(Y)\rightarrow\mathrm{q\text{-}Coh}(X) \]
by
\[ \mathcal{N}_{\pi}(\mathcal{G}):=\mathcal{N}_{\pi_1}(\mathcal{G}\mid_{X_1})\oplus\cdots\oplus\mathcal{N}_{\pi_n}(\mathcal{G}\mid_{X_n}) \]
for each $\mathcal{G}\in\mathrm{q\text{-}Coh}(Y)$.
\end{defn}

By the arguments so far, we obtain the following:

\begin{prop}\label{GeneralthetaProp}
Let $\pi:Y\rightarrow X$ be a finite \'{e}tale covering.

{\rm (i)} For any morphism $f:X^{\prime}\rightarrow X$, if we take the pull-back\[
\xy
(-10,8)*+{Y^{\prime}}="0";
(10,8)*+{Y}="2";
(-10,-8)*+{X^{\prime}}="4";
(10,-8)*+{X}="6";
(0,0)*+{\square}="8";
(13,-9)*+{,}="10";
{\ar^{g} "0";"2"};
{\ar_{\pi^{\prime}} "0";"4"};
{\ar_{f} "4";"6"};
{\ar^{\pi} "2";"6"};
\endxy
\]
then there exists a natural monoidal isomorphism
\[ \theta^f_{\pi}:f^{\ast}\circ\mathcal{N}_{\pi}\overset{\cong}{\longrightarrow}\mathcal{N}_{\pi^{\prime}}\circ g^{\ast}. \]

{\rm (ii)} For any other morphism $f^{\prime}:X^{\prime\prime}\rightarrow X^{\prime}$, if we take the pull-back
\[
\xy
(-10,8)*+{Y^{\prime\prime}}="0";
(10,8)*+{Y^{\prime}}="2";
(-10,-8)*+{X^{\prime\prime}}="4";
(10,-8)*+{X^{\prime}}="6";
(0,0)*+{\square}="8";
(13,-9)*+{,}="10";
{\ar^{g^{\prime}} "0";"2"};
{\ar_{\pi^{\prime\prime}} "0";"4"};
{\ar_{f^{\prime}} "4";"6"};
{\ar^{\pi^{\prime}} "2";"6"};
\endxy
\]
then we have
\[ \theta^{f\circ f^{\prime}}_{\pi}=(\theta^{f^{\prime}}_{\pi^{\prime}}\circ g^{\ast})\cdot(f^{\prime\ast}\circ\theta^f_{\pi}). \]
\end{prop}
\begin{pf}
This immediately follows from Proposition \ref{ConstthetaProp} and Corollary \ref{ConstthetaCor}.
\end{pf}

\begin{rem}
$\mathcal{N}_{\pi}$ is uniquely determined up to a natural monoidal isomorphism, by Definition \ref{TrivNormDef} and Proposition \ref{GeneralthetaProp}.
\end{rem}

\begin{prop}\label{Prop3.21}
Let $\pi:Y\rightarrow X$ be a finite \'{e}tale covering of constant degree $d$.
For any positive integer $m$, we have $:$
\[ \mathcal{N}_{\pi}(\mathcal{O}_Y^{\oplus m})\cong\mathcal{O}_X^{\oplus m^d} \]
\end{prop}
\begin{pf}
Take an fpqc morphism $f:X^{\prime}\rightarrow X$ trivializing $\pi$:
\[
\xy
(-18,6.2)*+{\underset{d}{\textstyle{\coprod}}X^{\prime}=}="-1";
(-8,8)*+{Y^{\prime}}="0";
(8,8)*+{Y}="2";
(-8,-8)*+{X^{\prime}}="4";
(8,-8)*+{X}="6";
(0,0)*+{\square}="8";
{\ar^{g} "0";"2"};
{\ar_{\nabla=\pi^{\prime}} "0";"4"};
{\ar_{f} "4";"6"};
{\ar^{\pi} "2";"6"};
\endxy
\]
Then we have an isomorphism
\begin{align*}
\beta:\mathcal{N}_{\pi^{\prime}}g^{\ast}(\mathcal{O}_Y^{\oplus m})&\cong\mathcal{N}_{\pi^{\prime}}(\mathcal{O}_{Y^{\prime}}^{\oplus m})\\
&=(\mathcal{O}_{Y^{\prime}}^{\oplus m})\mid_{X^{\prime}_1}\underset{X^{\prime}}{\otimes}\cdots\underset{X^{\prime}}{\otimes}(\mathcal{O}_{Y^{\prime}}^{\oplus m})\mid_{X^{\prime}_d}\\
&\cong\mathcal{O}_{X^{\prime}}^{\oplus m}\underset{X^{\prime}}{\otimes}\cdots\underset{X^{\prime}}{\otimes}\mathcal{O}_{X^{\prime}}^{\oplus m}\\
&\cong\mathcal{O}_{X^{\prime}}^{\oplus m^d}\cong f^{\ast}(\mathcal{O}_{X^{\prime}}^{\oplus m^d}).
\end{align*}
This $\beta$ satisfies the commutativity of
\[
\xy
(-32,8)*+{p_1^{\ast}\mathcal{N}_{\pi^{\prime}}g^{\ast}(\mathcal{O}_Y^{\oplus m})}="0";
(32,8)*+{p_1^{\ast}f^{\ast}(\mathcal{O}_{X^{\prime}}^{\oplus m^d})}="2";
(-32,-8)*+{p_2^{\ast}\mathcal{N}_{\pi^{\prime}}g^{\ast}(\mathcal{O}_Y^{\oplus m})}="4";
(32,-8)*+{p_2^{\ast}f^{\ast}(\mathcal{O}_{X^{\prime}}^{\oplus m^d})}="6";
{\ar^{p_1^{\ast}\beta} "0";"2"};
{\ar_{\phi}^{\cong} "0";"4"};
{\ar_{p_2^{\ast}\beta} "4";"6"};
{\ar^{\mathrm{can.}}_{\cong} "2";"6"};
{\ar@{}|\circlearrowright "0";"6"};
\endxy
\]
where $\phi:=\phi_{\mathcal{O}_Y^{\oplus m}}$ is the isomorphism defined in Claim \ref{Claim11}.
\end{pf}

\begin{cor}\label{Azumaya-1}
Let $\pi:Y\rightarrow X$ be a finite \'{e}tale covering.
If $\mathcal{E}\in\mathrm{q\text{-}Coh}(Y)$ is locally free of finite rank, then so is $\mathcal{N}_{\pi}(\mathcal{E})\in\mathrm{q\text{-}Coh}(X)$.
\end{cor}
\begin{pf}
By Proposition \ref{GeneralthetaProp}, we may assume $X$ is affine and connected.
Then $Y$ is also affine, and $\pi$ is of constant degree. Remark $\mathcal{E}$ is locally free of finite rank if and only if there is an integer $m$ and an epimorphism
\[ s:\mathcal{O}_Y^{\oplus m}\rightarrow \mathcal{E}. \]
Take an fpqc morphism $f:X^{\prime}\rightarrow X$ trivializing $\pi$ :
\[
\xy
(-18,6.2)*+{\underset{d}{\textstyle{\coprod}}X^{\prime}=}="-1";
(-8,8)*+{Y^{\prime}}="0";
(8,8)*+{Y}="2";
(-8,-8)*+{X^{\prime}}="4";
(8,-8)*+{X}="6";
(0,0)*+{\square}="8";
{\ar^{g} "0";"2"};
{\ar_{\nabla=\pi^{\prime}} "0";"4"};
{\ar_{f} "4";"6"};
{\ar^{\pi} "2";"6"};
\endxy
\]

By the definition of $\mathcal{N}_{\nabla}$, it can be easily seen that $\mathcal{N}_{\nabla}g^{\ast}(s)$ becomes epimorphic. Thus $f^{\ast}\mathcal{N}_{\pi}(s)$ is epimorphic.

Since $f$ is fully faithful, $\mathcal{N}_{\pi}(s):\mathcal{N}_{\pi}(\mathcal{O}_Y^{\oplus m})\rightarrow \mathcal{N}_{\pi}(\mathcal{E})$ also becomes epimorphic.
\end{pf}

\section{Norm maps}

\subsection{Norm map for the Brauer group}

\begin{defn}\label{deltapiDef}
Let $\pi:Y\rightarrow X$ be a finite \'{e}tale covering.
For any $\mathcal{F},\mathcal{G}\in\mathrm{q\text{-}Coh}(Y)$, we define a morphism
\[ \delta_{\pi}=\delta_{\pi,(\mathcal{F},\mathcal{G})}:\mathcal{N}_{\pi}\mathcal{H}om_{\mathcal{O}_Y}(\mathcal{F},\mathcal{G})\rightarrow\mathcal{H}om_{\mathcal{O}_X}(\mathcal{N}_{\pi}\mathcal{F},\mathcal{N}_{\pi}\mathcal{G}) \]
as follows $:$

Let
\[ e=\mathrm{ev}_{\mathcal{F},\mathcal{G}}:\mathcal{H}om_{\mathcal{O}_Y}(\mathcal{F},\mathcal{G})\underset{\mathcal{O}_Y}{\otimes}\mathcal{F}\rightarrow \mathcal{G} \]
be the evaluation morphism, i.e., the morphism corresponding to $\mathrm{id}_{\mathcal{H}om_{\mathcal{O}_Y}(\mathcal{F},\mathcal{G})}$ under the adjoint isomorphism.

Define $\xi_{\pi}$ as the composition
\[ \xi_{\pi}:=(\mathcal{N}_{\pi}\mathcal{H}om_{\mathcal{O}_Y}(\mathcal{F},\mathcal{G})\underset{\mathcal{O}_Y}{\otimes}\mathcal{N}_{\pi}\mathcal{F}\overset{\cong}{\longrightarrow}\mathcal{N}_{\pi}(\mathcal{H}om_{\mathcal{O}_Y}(\mathcal{F},\mathcal{G})\underset{\mathcal{O}_Y}{\otimes}\mathcal{F})\overset{\mathcal{N}_{\pi}(e)}{\longrightarrow}\mathcal{N}_{\pi}(\mathcal{G})). \]
By the adjoint isomorphism
\begin{align*}
&\mathrm{Hom}_{\mathcal{O}_X}(\mathcal{N}_{\pi}\mathcal{H}om_{\mathcal{O}_Y}(\mathcal{F},\mathcal{G})\underset{\mathcal{O}_X}{\otimes}\mathcal{N}_{\pi}\mathcal{F},\mathcal{N}_{\pi}\mathcal{G})\\
&\qquad\qquad\overset{\cong}{\longrightarrow}\mathrm{Hom}_{\mathcal{O}_X}(\mathcal{N}_{\pi}\mathcal{H}om_{\mathcal{O}_Y}(\mathcal{F},\mathcal{G}),\mathcal{H}om_{\mathcal{O}_X}(\mathcal{N}_{\pi}\mathcal{F},\mathcal{N}_{\pi}\mathcal{G})),
\end{align*}
we obtain $\delta_{\pi}$ corresponding to $\xi_{\pi}$.
\end{defn}

\begin{rem}
To define $\delta_{\pi}$, we only used the monoidality of $\mathcal{N}_{\pi}$.
In fact, for any monoidal functor $F:\mathcal{C}\rightarrow\mathcal{D}$ between closed symmetric monoidal categories, we can define a natural transformation
\[ \delta_{F}:F[-,-]_{\mathcal{C}}\rightarrow[F(-),F(-)]_{\mathcal{D}}, \]
where $[-,-]_{\mathcal{C}}$ and $[-,-]_{\mathcal{D}}$ are the right adjoint of $\otimes_{\mathcal{C}}$ and $\otimes_{\mathcal{D}}$, respectively.
\end{rem}

The following proposition also follows from general arguments on monoidal functors. We omit its proof.
\begin{prop}
In Definition \ref{deltapiDef}, if $\mathcal{F}=\mathcal{G}$, then
\[ \delta_{\pi}:\mathcal{N}_{\pi}\mathcal{E}nd_{\mathcal{O}_Y}(\mathcal{F})\rightarrow\mathcal{E}nd_{\mathcal{O}_X}(\mathcal{N}_{\pi}\mathcal{F}) \]
is a monoid morphism.
\end{prop}

\begin{rem}\label{IthetaRem}
Let $\pi:Y\rightarrow X$ be a finite \'{e}tale covering, and take the pull-back by a morphism $f:X^{\prime}\rightarrow X$ $:$
\[
\xy
(-10,8)*+{Y^{\prime}}="0";
(10,8)*+{Y}="2";
(-10,-8)*+{X^{\prime}}="4";
(10,-8)*+{X}="6";
(0,0)*+{\square}="8";
{\ar^{g} "0";"2"};
{\ar_{\pi^{\prime}} "0";"4"};
{\ar_{f} "4";"6"};
{\ar^{\pi} "2";"6"};
\endxy
\]

Let $\mathcal{F},\mathcal{G}\in\mathrm{q\text{-}Coh}(Y)$.
From $\theta^f_{\pi}$, we obtain an isomorphism
\[ I_{\theta}:\mathcal{H}om_{\mathcal{O}_{X^{\prime}}}(\mathcal{N}_{\pi^{\prime}}g^{\ast}\mathcal{F},\mathcal{N}_{\pi^{\prime}}g^{\ast}\mathcal{G})\overset{\cong}{\longrightarrow}\mathcal{H}om_{\mathcal{O}_{X^{\prime}}}(f^{\ast}\mathcal{N}_{\pi}\mathcal{F},f^{\ast}\mathcal{N}_{\pi}\mathcal{G}) \]
such that for any $\mathcal{E}\in\mathrm{q\text{-}Coh}(X^{\prime})$, the following diagram is commutative $:$
\[
\xy
(-36,18)*+{\mathrm{Hom}_{\mathcal{O}_{X^{\prime}}}(\mathcal{E},\mathcal{H}om(\mathcal{N}_{\pi^{\prime}}g^{\ast}\mathcal{F},\mathcal{N}_{\pi^{\prime}}g^{\ast}\mathcal{G}))}="1";
(36,18)*+{\mathrm{Hom}_{\mathcal{O}_{X^{\prime}}}(\mathcal{E},\mathcal{H}om(f^{\ast}\mathcal{N}_{\pi}\mathcal{F},f^{\ast}\mathcal{N}_{\pi}\mathcal{G}))}="2";
(-42,-3)*+{\mathrm{Hom}_{\mathcal{O}_{X^{\prime}}}(\mathcal{E}\otimes\mathcal{N}_{\pi^{\prime}}g^{\ast}\mathcal{F},\mathcal{N}_{\pi^{\prime}}g^{\ast}\mathcal{G})}="3";
(42,-3)*+{\mathrm{Hom}_{\mathcal{O}_{X^{\prime}}}(\mathcal{E}\otimes f^{\ast}\mathcal{N}_{\pi}\mathcal{F},f^{\ast}\mathcal{N}_{\pi}\mathcal{G})}="4";
(0,-20)*+{\mathrm{Hom}_{\mathcal{O}_{X^{\prime}}}(\mathcal{E}\otimes f^{\ast}\mathcal{N}_{\pi}\mathcal{F},\mathcal{N}_{\pi^{\prime}}g^{\ast}\mathcal{G})}="5";
(0,22)*+{}="6";
{\ar^{I_{\theta}\circ -} "1";"2"};
{\ar^{\cong}_{\mathrm{adj.}} "1";"3"};
{\ar_{\cong}^{\mathrm{adj.}} "2";"4"};
{\ar_{-\circ(\mathrm{id}\otimes\theta^f_{\pi})} "3";"5"};
{\ar^{\theta^f_{\pi}\circ -} "4";"5"};
{\ar@{}|\circlearrowright "5";"6"};
\endxy
\]
\end{rem}

\begin{prop}\label{IthetaProp}
In the notation of Remark \ref{IthetaRem}, assume $f$ is flat and $\mathcal{F}$ is locally free of finite rank. Remark there exist canonical natural isomorphisms
\begin{align*}
c_1&:g^{\ast}\mathcal{H}om_{\mathcal{O}_Y}(\mathcal{F},\mathcal{G})\overset{\cong}{\longrightarrow}\mathcal{H}om_{\mathcal{O}_{Y^{\prime}}}(g^{\ast}\mathcal{F},g^{\ast}\mathcal{G}),\\
c_2&:f^{\ast}\mathcal{H}om_{\mathcal{O}_X}(\mathcal{N}_{\pi}\mathcal{F},\mathcal{N}_{\pi}\mathcal{G})\overset{\cong}{\longrightarrow}\mathcal{H}om_{\mathcal{O}_{X^{\prime}}}(f^{\ast}\mathcal{N}_{\pi}\mathcal{F},f^{\ast}\mathcal{N}_{\pi}\mathcal{G}).
\end{align*}

Then, the following diagram is commutative $:$
\[
\xy
(-30,16)*+{f^{\ast}\mathcal{N}_{\pi}\mathcal{H}om_{\mathcal{O}_Y}(\mathcal{F},\mathcal{G})}="10";
(30,16)*+{f^{\ast}\mathcal{H}om_{\mathcal{O}_X}(\mathcal{N}_{\pi}\mathcal{F},\mathcal{N}_{\pi}\mathcal{G})}="12";
(-30,0)*+{\mathcal{N}_{\pi^{\prime}}g^{\ast}\mathcal{H}om_{\mathcal{O}_Y}(\mathcal{F},\mathcal{G})}="20";
(30,0)*+{\mathcal{H}om_{\mathcal{O}_{X^{\prime}}}(f^{\ast}\mathcal{N}_{\pi}\mathcal{F},f^{\ast}\mathcal{N}_{\pi}\mathcal{G})}="22";
(-30,-16)*+{\mathcal{N}_{\pi^{\prime}}\mathcal{H}om_{\mathcal{O}_{Y^{\prime}}}(g^{\ast}\mathcal{F},g^{\ast}\mathcal{G})}="30";
(30,-16)*+{\mathcal{H}om_{\mathcal{O}_{X^{\prime}}}(\mathcal{N}_{\pi^{\prime}}g^{\ast}\mathcal{F},\mathcal{N}_{\pi^{\prime}}g^{\ast}\mathcal{G})}="32";
{\ar^{f^{\ast}\delta_{\pi}} "10";"12"};
{\ar_{\theta^f_{\pi}}^{\cong} "10";"20"};
{\ar^{c_2}_{\cong} "12";"22"};
{\ar_{\mathcal{N}_{\pi^{\prime}}(c_1)}^{\cong} "20";"30"};
{\ar_{I_{\theta}}^{\cong} "32";"22"};
{\ar_{\delta_{\pi^{\prime}}} "30";"32"};
{\ar@{}|\circlearrowright "10";"32"};
\endxy
\]
\end{prop}
\begin{pf}
Put
\begin{align*}
e:=\mathrm{ev}_{\mathcal{F},\mathcal{G}}&:\mathcal{H}om_{\mathcal{O}_Y}(\mathcal{F},\mathcal{G})\otimes_{\mathcal{O}_Y}\mathcal{F}\rightarrow\mathcal{G},\\
e^{\prime}:=\mathrm{ev}_{g^{\ast}\mathcal{F},g^{\ast}\mathcal{G}}&:\mathcal{H}om_{\mathcal{O}_{Y^{\prime}}}(g^{\ast}\mathcal{F},g^{\ast}\mathcal{G})\otimes_{\mathcal{O}_{Y^{\prime}}}g^{\ast}\mathcal{F}\rightarrow g^{\ast}\mathcal{G}.
\end{align*}
Remark that
\[
\xy
(-32,10)*+{\mathcal{H}om_{\mathcal{O}_{Y^{\prime}}}(g^{\ast}\mathcal{F},g^{\ast}\mathcal{G})\otimes_{\mathcal{O}_{Y^{\prime}}}g^{\ast}\mathcal{F}}="3";
(32,10)*+{g^{\ast}(\mathcal{H}om_{\mathcal{O}_Y}(\mathcal{F},\mathcal{G})\otimes_{\mathcal{O}_Y}\mathcal{F})}="4";
(0,-8)*+{g^{\ast}\mathcal{G}}="5";
(0,12)*+{}="6";
{\ar_{e^{\prime}} "3";"5"};
{\ar^{\cong} "3";"4"};
{\ar^{g^{\ast}e} "4";"5"};
{\ar@{}|\circlearrowright "5";"6"};
\endxy
\]
is commutative. Put
\begin{align*}
u&:=I_{\theta}\circ\delta_{\pi^{\prime}}\circ\mathcal{N}_{\pi^{\prime}}(c_1)\circ\theta^f_{\pi},\\
v&:=c_2\circ f^{\ast}\delta_{\pi},
\end{align*}
and let $\mu$ and $\nu$ be their images under the adjoint isomorphism
\begin{eqnarray*}
&\mathrm{Hom}_{\mathcal{O}_{X^{\prime}}}(f^{\ast}\mathcal{N}_{\pi}\mathcal{H}om_{\mathcal{O}_Y}(\mathcal{F},\mathcal{G}),\mathcal{H}om_{\mathcal{O}_{X^{\prime}}}(f^{\ast}\mathcal{N}_{\pi}\mathcal{F},f^{\ast}\mathcal{N}_{\pi}\mathcal{G}))\\
&\qquad\qquad\qquad\overset{\cong}{\longrightarrow}\mathrm{Hom}_{\mathcal{O}_{X^{\prime}}}(f^{\ast}\mathcal{N}_{\pi}\mathcal{H}om_{\mathcal{O}_Y}(\mathcal{F},\mathcal{G})\otimes_{\mathcal{O}_{X^{\prime}}} f^{\ast}\mathcal{N}_{\pi}\mathcal{F},f^{\ast}\mathcal{N}_{\pi}\mathcal{G}),
\end{eqnarray*}
respectively.
It suffices to show $\mu=\nu$.

Put
\[ u_0:=\delta_{\pi^{\prime}}\circ\mathcal{N}_{\pi^{\prime}}(c_1)\circ\theta^f_{\pi}, \]
and let $\mu_0$ be its image under
\begin{eqnarray*}
&\mathrm{Hom}_{\mathcal{O}_{X^{\prime}}}(f^{\ast}\mathcal{N}_{\pi}\mathcal{H}om_{\mathcal{O}_Y}(\mathcal{F},\mathcal{G}),\mathcal{H}om_{\mathcal{O}_{X^{\prime}}}(\mathcal{N}_{\pi^{\prime}}g^{\ast}\mathcal{F},\mathcal{N}_{\pi^{\prime}}g^{\ast}\mathcal{G}))\\
&\qquad\qquad\qquad\overset{\cong}{\longrightarrow}\mathrm{Hom}_{\mathcal{O}_{X^{\prime}}}(f^{\ast}\mathcal{N}_{\pi}\mathcal{H}om_{\mathcal{O}_Y}(\mathcal{F},\mathcal{G})\otimes_{\mathcal{O}_{X^{\prime}}}\mathcal{N}_{\pi^{\prime}}g^{\ast}\mathcal{F},\mathcal{N}_{\pi^{\prime}}g^{\ast}\mathcal{G})
\end{eqnarray*}

Since $\mu=I_{\theta}\circ\mu_0$, by Remark \ref{IthetaRem}, we have
\begin{equation}
\mu_0\circ(\mathrm{id}\otimes\theta^f_{\pi})=\theta^f_{\pi}\circ\mu.
\label{LEQ1}
\end{equation}
\[
\xy
(-24,8)*+{f^{\ast}\mathcal{N}_{\pi}\mathcal{H}om_{\mathcal{O}_Y}(\mathcal{F},\mathcal{G})\otimes_{\mathcal{O}_{X^{\prime}}}f^{\ast}\mathcal{N}_{\pi}\mathcal{F}}="0";
(24,8)*+{f^{\ast}\mathcal{N}_{\pi}\mathcal{G}}="2";
(-24,-8)*+{f^{\ast}\mathcal{N}_{\pi}\mathcal{H}om_{\mathcal{O}_Y}(\mathcal{F},\mathcal{G})\otimes_{\mathcal{O}_{X^{\prime}}}\mathcal{N}_{\pi^{\prime}}g^{\ast}\mathcal{F}}="4";
(24,-8)*+{\mathcal{N}_{\pi^{\prime}}g^{\ast}\mathcal{G}}="6";
{\ar^<<<<<<<<<{\mu} "0";"2"};
{\ar_{\mathrm{id}\otimes\theta^f_{\pi}} "0";"4"};
{\ar_>>>>>>>>>{\mu_0} "4";"6"};
{\ar^{\theta^f_{\pi}} "2";"6"};
{\ar@{}|\circlearrowright "0";"6"};
\endxy
\]

By the definition of $\xi_{\pi^{\prime}}$ and the naturality of the adjoint isomorphism, we can show easily
\begin{equation}
\mu_0=\xi_{\pi^{\prime}}\circ(\mathcal{N}_{\pi^{\prime}}(c_1)\otimes\mathrm{id})\circ(\theta^f_{\pi}\otimes\mathrm{id}).
\label{LEQ2}
\end{equation}

\[
\xy
(-36,8)*+{f^{\ast}\mathcal{N}_{\pi}\mathcal{H}om(\mathcal{F},\mathcal{G})\otimes\mathcal{N}_{\pi^{\prime}}g^{\ast}\mathcal{F}}="0";
(36,8)*+{\mathcal{N}_{\pi^{\prime}}g^{\ast}\mathcal{G}}="2";
(-36,-8)*+{\mathcal{N}_{\pi^{\prime}}g^{\ast}\mathcal{H}om(\mathcal{F},\mathcal{G})\otimes\mathcal{N}_{\pi^{\prime}}g^{\ast}\mathcal{F}}="4";
(36,-8)*+{\mathcal{N}_{\pi^{\prime}}\mathcal{H}om(g^{\ast}\mathcal{F},g^{\ast}\mathcal{G})\otimes\mathcal{N}_{\pi^{\prime}}g^{\ast}\mathcal{F}}="6";
{\ar^>>>>>>>>>>>>>>>>>>>>>>>>>>>>{\mu_0} "0";"2"};
{\ar_{\mathcal{N}_{\pi^{\prime}}(c_1)\otimes\mathrm{id}} "4";"6"};
{\ar_{\theta^f_{\pi}\otimes\mathrm{id}} "0";"4"};
{\ar_{\xi_{\pi^{\prime}}} "6";"2"};
{\ar@{}|\circlearrowright "0";"6"};
\endxy
\]

On the other hand, we have a commutative diagram
\begin{equation}
\xy
(-24,8)*+{f^{\ast}\mathcal{N}_{\pi}\mathcal{H}om(\mathcal{F},\mathcal{G})\otimes f^{\ast}\mathcal{N}_{\pi}\mathcal{F}}="0";
(0,-12)*+{f^{\ast}(\mathcal{N}_{\pi}\mathcal{H}om(\mathcal{F},\mathcal{G})\otimes\mathcal{N}_{\pi}\mathcal{F})}="2";
(25,-13)*+{.}="3";
(24,8)*+{f^{\ast}\mathcal{N}_{\pi}\mathcal{G}}="4";
(0,14)*+{}="6";
{\ar_{\cong} "0";"2"};
{\ar^>>>>>>>>>>>>>{\nu} "0";"4"};
{\ar_{f^{\ast}\xi_{\pi}} "2";"4"};
{\ar@{}|\circlearrowright "2";"6"};
\endxy
\label{LEQ3}
\end{equation}

Moreover, since $\theta^f_{\pi}$ is a natural monoidal transformation, the following diagram is commutative:
\begin{equation}
\xy
(-34,16)*+{f^{\ast}\mathcal{N}_{\pi}\mathcal{H}om(\mathcal{F},\mathcal{G})\otimes f^{\ast}\mathcal{N}_{\pi}\mathcal{F}}="10";
(34,16)*+{\mathcal{N}_{\pi^{\prime}}g^{\ast}\mathcal{H}om(\mathcal{F},\mathcal{G})\otimes\mathcal{N}_{\pi^{\prime}}g^{\ast}\mathcal{F}}="12";
(-34,0)*+{f^{\ast}\mathcal{N}_{\pi}(\mathcal{H}om(\mathcal{F},\mathcal{G})\otimes\mathcal{G})}="20";
(34,0)*+{\mathcal{N}_{\pi^{\prime}}g^{\ast}(\mathcal{H}om(\mathcal{F},\mathcal{G})\otimes\mathcal{G})}="22";
(-34,-16)*+{f^{\ast}\mathcal{N}_{\pi}\mathcal{G}}="30";
(34,-16)*+{\mathcal{N}_{\pi^{\prime}}g^{\ast}\mathcal{G}}="32";
{\ar^{\theta^f_{\pi}\otimes\theta^f_{\pi}} "10";"12"};
{\ar_{\theta^f_{\pi}} "20";"22"};
{\ar_{\theta^f_{\pi}} "30";"32"};
{\ar_{\cong} "10";"20"};
{\ar^{\cong} "12";"22"};
{\ar_{f^{\ast}\mathcal{N}_{\pi}(e)} "20";"30"};
{\ar^{\mathcal{N}_{\pi^{\prime}}g^{\ast}(e)} "22";"32"};
{\ar@{}|\circlearrowright "10";"22"};
{\ar@{}|\circlearrowright "20";"32"};
\endxy
\label{LEQ4}
\end{equation}
From $(\ref{LEQ1}),(\ref{LEQ2}),(\ref{LEQ3}),(\ref{LEQ4})$, we obtain $\mu=\nu$.
\end{pf}

\begin{cor}
Let $\pi:Y\rightarrow X$ be a finite \'{e}tale covering, and let $\mathcal{F},\mathcal{G}\in\mathrm{q\text{-}Coh}(Y)$. If $\mathcal{F}$ is locally free of finite rank, then
\[ \delta_{\pi}:\mathcal{N}_{\pi}\mathcal{H}om_{\mathcal{O}_Y}(\mathcal{F},\mathcal{G})\rightarrow\mathcal{H}om_{\mathcal{O}_X}(\mathcal{N}_{\pi}\mathcal{F},\mathcal{N}_{\pi}\mathcal{G}) \]
is an isomorphism.
\end{cor}
\begin{pf}
Let $U\subset X$ be any open subscheme.
Put $V:=\pi^{-1}(U)$ and let $\varpi:V\rightarrow U$ be the restriction of $\pi$ onto $V$.

By Proposition \ref{IthetaProp}, we have a commutative diagram
\[
\xy
(-36,8)*+{(\mathcal{N}_{\pi}\mathcal{H}om_{\mathcal{O}_Y}(\mathcal{F},\mathcal{G}))\mid_U}="0";
(36,8)*+{(\mathcal{H}om_{\mathcal{O}_X}(\mathcal{N}_{\pi}\mathcal{F},\mathcal{N}_{\pi}\mathcal{G}))\mid_U}="2";
(-36,-8)*+{\mathcal{N}_{\varpi}\mathcal{H}om_{\mathcal{O}_V}(\mathcal{F}\mid_V,\mathcal{G}\mid_V)}="4";
(36,-8)*+{\mathcal{H}om_{\mathcal{O}_U}(\mathcal{N}_{\varpi}(\mathcal{F}\mid_V),\mathcal{N}_{\varpi}(\mathcal{G}\mid_V))}="6";
{\ar^{(\delta_{\pi})\mid_U} "0";"2"};
{\ar_{\cong} "0";"4"};
{\ar_{\delta_{\varpi}} "4";"6"};
{\ar_{\cong} "2";"6"};
{\ar@{}|\circlearrowright "0";"6"};
\endxy
\]

Thus by taking an affine open cover of $X$, we may assume $X$ is affine and connected.

Moreover, again by Proposition \ref{IthetaProp}, replacing $X$ by its finite \'{e}tale covering $X^{\prime}\rightarrow X$, we may assume $Y$ is trivial over $X$, i.e.,
\begin{align*}
Y&=\underset{1\le k\le d}{\textstyle{\coprod}}X_k\quad(X_k=X\quad (1\le \forall k\le d)),\\
\pi&=\nabla_{X,d}.
\end{align*}
Since $Y$ is affine, any $\mathcal{F}\in\mathrm{q\text{-}Coh}(Y)$ can be identified with $\Gamma(Y,\mathcal{F})$, which is a $\Gamma(Y,\mathcal{O}_Y)$-module. Similarly for the sheaves on $X$.

Under this identification, for any $\mathcal{F},\mathcal{G}\in\mathrm{q\text{-}Coh}(Y)$, $\mathcal{H}om_{\mathcal{O}_Y}(\mathcal{F},\mathcal{G})$ is regarded as $\mathrm{Hom}_{\mathcal{O}_Y}(\mathcal{F},\mathcal{G})$ and
\[ e=\mathrm{ev}_{\mathcal{F},\mathcal{G}}:\mathcal{H}om_{\mathcal{O}_Y}(\mathcal{F},\mathcal{G})\otimes_{\mathcal{O}_Y}\mathcal{F}\rightarrow\mathcal{G} \]
is given by
\[ e(\varphi\otimes x)=\varphi(x)\quad (\forall\varphi\in\mathrm{Hom}_{\mathcal{O}_Y}(\mathcal{F},\mathcal{G}),\forall x\in\Gamma(Y,\mathcal{F})). \]
Similarly, it can be easily seen that $\xi_{\pi}$
\begin{align*}
&\qquad(\mathcal{H}om(\mathcal{F},\mathcal{G})\mid_{X_1}\underset{\mathcal{O}_X}{\otimes}\cdots\underset{\mathcal{O}_X}{\otimes}\mathcal{H}om(\mathcal{F},\mathcal{G})\mid_{X_d})\underset{\mathcal{O}_X}{\otimes}(\mathcal{F}\mid_{X_1}\underset{\mathcal{O}_X}{\otimes}\cdots\underset{\mathcal{O}_X}{\otimes}\mathcal{F}\mid_{X_d})\\
&\hspace{6cm}\rotatebox{90}{$\cong$}
\\
&(\mathcal{H}om(\mathcal{F}\mid_{X_1},\mathcal{G}\mid_{X_1})\underset{\mathcal{O}_X}{\otimes}\cdots\underset{\mathcal{O}_X}{\otimes}\mathcal{H}om(\mathcal{F}\mid_{X_d},\mathcal{G}\mid_{X_d})\underset{\mathcal{O}_X}{\otimes}(\mathcal{F}\mid_{X_1}\underset{\mathcal{O}_X}{\otimes}\cdots\underset{\mathcal{O}_X}{\otimes}\mathcal{F}\mid_{X_d})\\
&\hspace{6cm}\downarrow
\\
&\hspace{4.4cm}\mathcal{G}\mid_{X_1}\underset{\mathcal{O}_X}{\otimes}\cdots\underset{\mathcal{O}_X}{\otimes}\mathcal{G}\mid_{X_d}
\end{align*}
is given by
\[ \xi_{\pi}((\varphi_1\otimes\cdots\otimes\varphi_d)\otimes(x_1\otimes\cdots\otimes x_d))=\varphi_1(x_1)\otimes\cdots\otimes\varphi_d(x_d) \]
for any $\varphi_k\in\mathrm{Hom}_{\mathcal{O}_{X_k}}(\mathcal{F}\mid_{X_k},\mathcal{G}\mid_{X_k})$ and $x_k\in\Gamma(X_k,\mathcal{F}\mid_{X_k})$ $(1\le\forall k\le d)$.

Correspondingly, $\delta_{\pi}$ is given by
\[
\xy
(-40,8)*+{\mathrm{Hom}_{\mathcal{O}_X}(\mathcal{F}\mid_{X_1},\mathcal{G}\mid_{X_1})\underset{\mathcal{O}_X}{\otimes}\cdots\underset{\mathcal{O}_X}{\otimes}\mathrm{Hom}_{\mathcal{O}_X}(\mathcal{F}\mid_{X_d},\mathcal{G}\mid_{X_d})\ni\varphi_1\otimes\cdots\otimes\varphi_d}="0";
(-40,-8)*+{\mathrm{Hom}_{\mathcal{O}_X}(\mathcal{F}\mid_{X_1}\underset{\mathcal{O}_X}{\otimes}\cdots\underset{\mathcal{O}_X}{\otimes}\mathcal{F}\mid_{X_d},\mathcal{G}\mid_{X_1}\underset{\mathcal{O}_X}{\otimes}\cdots\underset{\mathcal{O}_X}{\otimes}\mathcal{G}\mid_{X_d})\ni\varphi_1\otimes\cdots\otimes\varphi_d,}="2";
{\ar (-52,4);(-52,-3)};
{\ar@{|->} (8,4);(8,-3)};
\endxy
\]
which is isomorphic.
\end{pf}

\begin{cor}\label{deltapiCor}
In particular, for any locally free $\mathcal{E}\in\mathrm{q\text{-}Coh}(Y)$ of finite rank,
\[ \delta_{\pi}:\mathcal{N}_{\pi}(\mathcal{E}nd_{\mathcal{O}_Y}(\mathcal{E}))\rightarrow\mathcal{E}nd_{\mathcal{O}_X}(\mathcal{N}_{\pi}(\mathcal{E})) \]
is an isomorphism of $\mathcal{O}_X$-algebras.
\end{cor}

\begin{lem}
\label{ec}
Let $\pi:Y\rightarrow X$ be a finite \'{e}tale covering. For any surjective \'{e}tale morphism $g:V\rightarrow Y$, there exists a surjective \'{e}tale morphism $f:U\rightarrow X$ such that
\[ \mathrm{pr}_Y:U\times_XY\rightarrow Y \]
factors through $g$.
\end{lem}
\begin{pf}
We may assume $\pi$ is  of constant degree $d$.
Let $x:X^{\prime}\rightarrow X$ be a surjective \'{e}tale covering which trivializes $\pi\, :$
\[
\xy
(-24,6)*+{\underset{1\le k\le d}{\textstyle{\coprod}}X^{\prime}_k=}="8";
(-10,8)*+{Y^{\prime}}="10";
(10,8)*+{Y}="12";
(-10,-8)*+{X^{\prime}}="20";
(10,-8)*+{X}="22";
(0,0)*+{\square}="1";
{\ar^{y} "10";"12"};
{\ar_{\nabla=\pi^{\prime}} "10";"20"};
{\ar_{x} "20";"22"};
{\ar^{\pi} "12";"22"};
\endxy
\]

Pulling $g$ back by $y$, we obtain
\[
\xy
(-20,8)*+{V^{\prime}_k\ }="10";
(0,8)*+{V^{\prime}}="12";
(-20,-8)*+{X^{\prime}_k\,}="20";
(0,-8)*+{Y^{\prime}}="22";
(20,8)*+{V}="14";
(20,-8)*+{Y}="24";
(-10,0)*+{\square}="1";
(10,0)*+{\square}="2";
{\ar@{^(->}^{\iota_k} "10";"12"};
{\ar_{g^{\prime}_k} "10";"20"};
{\ar@{^(->}_{\eta_k} "20";"22"};
{\ar^{g^{\prime}} "12";"22"};
{\ar^{v} "12";"14"};
{\ar_{y} "22";"24"};
{\ar_{g} "14";"24"};
\endxy
\]
where $\iota_k$ and $\eta_k$ $(1\le k\le d)$ are open immersions.

If we put $U:=V^{\prime}_1\times_{X^{\prime}}\cdots\times_{X^{\prime}}V^{\prime}_d$, we obtain a surjective \'{e}tale morphism $f:U\rightarrow X^{\prime}$ which makes the following diagram commutative for any $1\le k\le d\, :$
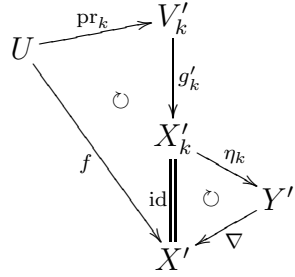
\[
\xy
(-10,12)*+{U}="0";
(10,16)*+{V^{\prime}_k}="2";
(10,0)*+{X^{\prime}_k}="4";
(10,-16)*+{X^{\prime}}="6";
(24,-8)*+{Y^{\prime}}="8";
(-4,10)*+{}="10";
(20,-16)*+{}="12";
{\ar^{\mathrm{pr}_k} "0";"2"};
{\ar_{f} "0";"6"};
{\ar^{g^{\prime}_k} "2";"4"};
{\ar@{=}_{\mathrm{id}} "4";"6"};
{\ar^{\eta_k} "4";"8"};
{\ar^{\nabla} "8";"6"};
{\ar@{}|\circlearrowright "4";"10"};
{\ar@{}|\circlearrowright "4";"12"};
\endxy
\]
Pulling back $f$, we obtain:
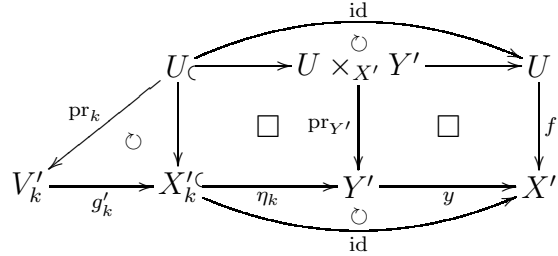
\[
\xy
(-20,8)*+{U }="10";
(4,8)*+{U\times_{X^{\prime}}Y^{\prime}}="12";
(4,14)*+{}="13";
(-20,-8)*+{X^{\prime}_k}="20";
(-32,4)*+{}="21";
(-40,-8)*+{V^{\prime}_k}="0";
(4,-8)*+{Y^{\prime}}="22";
(4,-16)*+{}="23";
(28,8)*+{U}="14";
(28,-8)*+{X^{\prime}}="24";
(-8,0)*+{\square}="1";
(16,0)*+{\square}="2";
{\ar@{_(->} "10";"12"};
{\ar^{} "10";"20"};
{\ar@{^(->}_{\eta_k} "20";"22"};
{\ar_{\mathrm{pr}_{Y^{\prime}}} "12";"22"};
{\ar_{\mathrm{pr}_k} "10";"0"};
{\ar_{g^{\prime}_k} "0";"20"};
{\ar "12";"14"};
{\ar_{y} "22";"24"};
{\ar^{f} "14";"24"};
{\ar@{}|\circlearrowright "20";"21"};
{\ar@{}|\circlearrowright "12";"13"};
{\ar@{}|\circlearrowright "22";"23"};
{\ar@/^1.40pc/^{\mathrm{id}} "10";"14"};
{\ar@/_1.40pc/_{\mathrm{id}} "20";"24"};
\endxy
\]

Thus, if we put
\[ p:=\underset{k}{\textstyle{\coprod}}\,\mathrm{pr}_k:U\times_{X^{\prime}}Y^{\prime}\rightarrow V^{\prime}, \]
then we have the following commutative diagram :
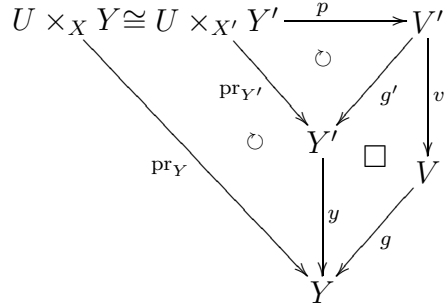
\[
\xy
(-14,6)*+{U\times_{X^{\prime}}Y^{\prime}}="0";
(-18,10)*+{}="1";
(14,6)*+{V^{\prime}}="2";
(0,-10)*+{Y^{\prime}}="4";
(0,12)*+{}="5";
(-34,6)*+{U\times_XY}="10";
(-25,6)*+{\cong}="11";
(14,-14)*+{V}="12";
(0,-30)*+{Y}="14";
(7,-12)*+{\square}="20";
{\ar^{p} "0";"2"};
{\ar^{g^{\prime}} "2";"4"};
{\ar_{\mathrm{pr}_{Y^{\prime}}} "0";"4"};
{\ar_{\mathrm{pr}_{Y}} "10";"14"};
{\ar^{y} "4";"14"};
{\ar^{g} "12";"14"};
{\ar@{}|\circlearrowright "1";"14"};
{\ar@{}|\circlearrowright "4";"5"};
{\ar^{v} "2";"12"};
\endxy
\]
This is what we wanted to show.
\end{pf}

\begin{prop}
\label{Azumaya-2}
Let $\pi:Y\rightarrow X$ be a finite \'{e}tale covering.
If $\mathcal{A}\in\mathrm{q\text{-}Coh}(Y)$ is an Azumaya algebra on $Y$, then $\mathcal{N}_{\pi}(\mathcal{A})$ becomes an Azumaya algebra on $X$.
\end{prop}
\begin{pf}
Since $\mathcal{N}_{\pi}$ is monoidal, $\mathcal{N}_{\pi}(\mathcal{A})$ becomes an $\mathcal{O}_X$-algebra. By Proposition \ref{Azumaya-1}, $\mathcal{N}_{\pi}(\mathcal{A})$ is locally free of finite rank.
Let $g:V\rightarrow Y$ be a surjective \'{e}tale morphism such that there exists a locally free sheaf $\mathcal{E}\in\mathrm{q\text{-}Coh}(Y)$ of finite rank, with an isomorphism of $\mathcal{O}_V$-algebras
\begin{equation}
g^{\ast}\mathcal{A}\cong g^{\ast}\mathcal{E}nd_{\mathcal{O}_Y}(\mathcal{E}).
\label{AzEq}
\end{equation}
By Lemma \ref{ec}, replacing $g$ if necessary, we may assume there exists a surjective \'{e}tale morphism $f:U\rightarrow X$ such that $g$ is the pull-back of $f$ by $\pi$:
\[
\xy
(-10,8)*+{V}="0";
(10,8)*+{Y}="2";
(-10,-8)*+{U}="4";
(10,-8)*+{X}="6";
(0,0)*+{\square}="8";
{\ar^{g} "0";"2"};
{\ar_{\varpi} "0";"4"};
{\ar_{f} "4";"6"};
{\ar^{\pi} "2";"6"};
\endxy
\]

By $(\ref{AzEq})$ and Corollary \ref{deltapiCor}, we obtain an isomorphism of $\mathcal{O}_U$-algebras
\begin{align*}
f^{\ast}\mathcal{N}_{\pi}(\mathcal{A})&\overset{\cong}{\rightarrow}\mathcal{N}_{\varpi}g^{\ast}\mathcal{A}\overset{\cong}{\rightarrow}\mathcal{N}_{\varpi}g^{\ast}(\mathcal{E}nd_{\mathcal{O}_Y}(\mathcal{E}))\\
&\overset{\cong}{\rightarrow}f^{\ast}\mathcal{N}_{\pi}(\mathcal{E}nd_{\mathcal{O}_Y}(\mathcal{E}))\overset{\cong}{\rightarrow}f^{\ast}(\mathcal{E}nd_{\mathcal{O}_X}(\mathcal{N}_{\pi}(\mathcal{E}))),
\end{align*}
which shows $\mathcal{N}_{\pi}(\mathcal{A})$ is an Azumaya algebra on $X$.
\end{pf}

\begin{cor}\label{MainResult1}
Let $\pi:Y\rightarrow X$ be a finite \'{e}tale covering. Norm functor $\mathcal{N}_{\pi}:\mathrm{q\text{-}Coh}(Y)\rightarrow\mathrm{q\text{-}Coh}(X)$ induces a group homomorphism
\[ N_{\pi}:\mathrm{Br}(Y)\rightarrow\mathrm{Br}(X), \]
which we call the {\it norm map}.
\end{cor}
\begin{pf}
This follows from Corollary \ref{deltapiCor} and Proposition \ref{Azumaya-2}.
\end{pf}

\subsection{Norm map for the cohomology group}

\begin{rem}
Remark there is a natural isomorphism
\[ \gamma_X:\Gamma(X,\mathbb{G}_{m,X})\overset{\cong}{\longrightarrow}\mathrm{Aut}_{\mathcal{O}_X}(\mathcal{O}_X) \]
for each scheme $X$. If $\pi:Y\rightarrow X$ is a finite \'{e}tale covering, from the norm functor $\mathcal{N}_{\pi}:\mathrm{q\text{-}Coh}(Y)\rightarrow\mathrm{q\text{-}Coh}(X)$, we obtain a group homomorphism
\[ %
\overline{\mathcal{N}_{\pi}}:%
\mathrm{Aut}_{\mathcal{O}_Y}(\mathcal{O}_Y)%
\overset{\mathcal{N}_{\pi}}{\longrightarrow}%
\mathrm{Aut}_{\mathcal{O}_X}(\mathcal{N}_{\pi}(\mathcal{O}_Y))%
\overset{\cong}{\longrightarrow}%
\mathrm{Aut}_{\mathcal{O}_X}(\mathcal{O}_X). \]
Thus we can define a group homomorphism
\[ N_{\pi}(X):=\gamma_X^{-1}\circ\overline{\mathcal{N}_{\pi}}\circ\gamma_Y:\Gamma(Y,\mathbb{G}_{m,Y})\rightarrow\Gamma(X,\mathbb{G}_{m,X}). \]
\end{rem}

\begin{prop}\label{NormHomProp}
Let $\pi:Y\rightarrow X$ be a finite \'{e}tale covering, and $f:U\rightarrow X$ be any \'{e}tale morphism of finite type.
Take the pull-back diagram
\[
\xy
(-8,8)*+{V}="0";
(8,8)*+{Y}="2";
(-8,-8)*+{U}="4";
(8,-8)*+{X}="6";
(11,-9)*+{.}="7";
(0,0)*+{\square}="9";
{\ar_{\varpi} "0";"4"};
{\ar^{\pi} "2";"6"};
{\ar^{g} "0";"2"};
{\ar_{f} "4";"6"};
\endxy
\]
We define
\[ N_{\pi}(U):\Gamma_{\mathrm{et}}(U,\pi_{\ast}\mathbb{G}_{m,Y})\rightarrow\Gamma_{\mathrm{et}}(U,\mathbb{G}_{m,X}) \]
by
\[ N_{\pi}(U):=N_{\varpi}(U):\Gamma(V,\mathcal{O}_V^{\times})\rightarrow\Gamma(U,\mathcal{O}_U^{\times}).\]
Then the set of group homomorphisms
\[ \{ N_{\pi}(U)\mid (U\overset{f}{\longrightarrow}X)\in X_{\mathrm{et}}\} \]
gives a homomorphism of abelian sheaves on $X_{\mathrm{et}}\, :$
\[ N_{\pi}:\pi_{\ast}\mathbb{G}_{m,Y}\rightarrow \mathbb{G}_{m,X} \]
\end{prop}
\begin{pf}
Let $f^{\prime}:U^{\prime}\rightarrow X$ be another \'{e}tale morphism of finite type, and $u:U\rightarrow U^{\prime}$ be an \'{e}tale morphism over $X$.
Take the pull-backs:
\[
\xy
(-16,8)*+{V^{\prime}}="0";
(0,8)*+{V}="2";
(0,16)*+{}="3";
(16,8)*+{Y}="4";
(-16,-8)*+{U^{\prime}}="10";
(0,-8)*+{U}="12";
(0,-16)*+{}="13";
(16,-8)*+{X}="14";
(-8,0)*+{\square}="30";
(8,0)*+{\square}="32";
{\ar^{v} "0";"2"};
{\ar^{g} "2";"4"};
{\ar_{u} "10";"12"};
{\ar_{f} "12";"14"};
{\ar^{\pi} "4";"14"};
{\ar^{\varpi} "2";"12"};
{\ar_{\varpi^{\prime}} "0";"10"};
{\ar@/^1.60pc/^{g^{\prime}} "0";"4"};
{\ar@/_1.60pc/_{f^{\prime}} "10";"14"};
{\ar@{}|\circlearrowright "2";"3"};
{\ar@{}|\circlearrowright "12";"13"};
\endxy
\]
It suffices to show the commutativity of
\[
\xy
(-20,8)*+{\Gamma_{\mathrm{et}}(V,\mathbb{G}_{m,Y})}="0";
(20,8)*+{\Gamma_{\mathrm{et}}(U,\mathbb{G}_{m,X})}="2";
(-20,-8)*+{\Gamma_{\mathrm{et}}(V^{\prime},\mathbb{G}_{m,Y})}="4";
(20,-8)*+{\Gamma_{\mathrm{et}}(U^{\prime},\mathbb{G}_{m,X})}="6";
(34,-9)*+{.}="7";
{\ar_{v^{\ast}} "0";"4"};
{\ar^{u^{\ast}} "2";"6"};
{\ar^{N_{\varpi}(U)} "0";"2"};
{\ar_{N_{\varpi^{\prime}}(U^{\prime})} "4";"6"};
{\ar@{}|\circlearrowright "0";"6"};
\endxy
\]
This immediately follows from the fact that $\theta$ is a natural monoidal isomorphism;
\[
\xy
(-32,0)*+{\mathrm{Aut}_{\mathcal{O}_V}(\mathcal{O}_V)}="0";
(-24,12)*+{\mathrm{Aut}_{\mathcal{O}_U}(\mathcal{N}_{\varpi}\mathcal{O}_V)}="2";(16,8)*+{\mathrm{Aut}_{\mathcal{O}_{U^{\prime}}}(u^{\ast}\mathcal{N}_{\varpi}\mathcal{O}_V)}="4";
(-24,-12)*+{\mathrm{Aut}_{\mathcal{O}_{V^{\prime}}}(v^{\ast}\mathcal{O}_V)}="6";(20,0)*+{}="5";
(16,-8)*+{\mathrm{Aut}_{\mathcal{O}_{U^{\prime}}}(\mathcal{N}_{\varpi^{\prime}}v^{\ast}\mathcal{O}_V)}="8";
(52,0)*+{\mathrm{Aut}_{\mathcal{O}_{U^{\prime}}}(\mathcal{O}_{U^{\prime}})}="10";
(8,0)*+{}="11";
{\ar_{v^{\ast}} "0";"6"};
{\ar^{u^{\ast}} "2";"4"};
{\ar^{\mathcal{N}_{\varpi}} "0";"2"};
{\ar_{\mathcal{N}_{\varpi^{\prime}}} "6";"8"};
{\ar_{c(\theta^u_{\varpi})} "4";"8"};
{\ar^{\cong} "4";"10"};
{\ar_{\cong} "8";"10"};
{\ar@{}|\circlearrowright "0";"5"};
{\ar@{}|\circlearrowright "10";"11"};
\endxy
\]
where $c(\theta^u_{\varpi})$ is the conjugation by $\theta^u_{\varpi}$.
\end{pf}

\begin{defn}
Let $\pi:Y\rightarrow X$ be a finite \'{e}tale covering.
By Proposition \ref{NormHomProp}, we obtain a homomorphism
\[ H^2_{\mathrm{et}}(N_{\pi}):H^2_{\mathrm{et}}(X,\pi_{\ast}\mathbb{G}_{m,Y})\rightarrow H^2_{\mathrm{et}}(X,\mathbb{G}_{m,X}). \]
We define the norm map for cohomology, as the composition of this map with the canonical isomorphism
\[ \mathfrak{c}^{-1}:H^2_{\mathrm{et}}(Y,\mathbb{G}_{m,Y})\overset{\cong}{\longrightarrow}H^2_{\mathrm{et}}(X,\pi_{\ast}\mathbb{G}_{m,Y}), \]
and abbreviately denote it by $N_{\pi}\,:$
\[ N_{\pi}:H^2_{\mathrm{et}}(Y,\mathbb{G}_{m,Y})\rightarrow H^2_{\mathrm{et}}(X,\mathbb{G}_{m,X}) \]
\end{defn}

\section{Compatibility of the norm maps}

In the following, we often assume that a scheme $X$ satisfies the following assumption:
\begin{assum}
\label{Assumption}
For any finite subset $F$ of $X$, there exists an affine open subscheme $U\subset X$ containing $F$.
\end{assum}

Remark that if $X$ satisfies Assumption \ref{Assumption}, then so does any finite \'{e}tale covering $Y$ over $X$.

\begin{rem}
Assumption \ref{Assumption} is only used in the proof of the next theorem. So if one can show it by another way, any of the succeeding results does not require Assumption \ref{Assumption}.
\end{rem}

\begin{thm}\label{TheoremToShow+}
For any finite \'{e}tale covering $\pi:Y\rightarrow X$, we have a commutative diagram
\[
\xy
(-18,8)*+{\mathrm{Br}(Y)}="0";
(18,8)*+{\mathrm{Br}(X)}="2";
(-18,4)*+{}="1";
(18,4)*+{}="3";
(-18,-8)*+{H^2_{\mathrm{et}}(Y,\mathbb{G}_{m,Y})}="4";
(18,-8)*+{H^2_{\mathrm{et}}(X,\mathbb{G}_{m,X}).}="6";
{\ar@{^(->}^{\chi_Y} "1";"4"};
{\ar@{^(->}^{\chi_X} "3";"6"};
{\ar^{N_{\pi}} "0";"2"};
{\ar_{N_{\pi}} "4";"6"};
{\ar@{}|\circlearrowright "0";"6"};
\endxy
\]
\end{thm}
\begin{pf}
By Assumption \ref{Assumption}, it suffices to show for the \v{C}ech cohomology. First, we briefly recall the construction of
\[ \chi_Y:\mathrm{Br}(Y)\hookrightarrow H^2_{\mathrm{et}}(Y,\mathbb{G}_{m,Y}) \]
using \v{C}ech cohomology (cf. \cite{Milne}). For any Azumaya algebra $\mathcal{A}$ on $Y$, there exists a surjective \'{e}tale morphism $g:V\rightarrow Y$, a locally free $\mathcal{E}\in \mathrm{q\text{-}Coh}(Y)$ of finite rank, and an isomorphism of $\mathcal{O}_V$-algebras
\[ \phi:g^{\ast}\mathcal{A}\overset{\cong}{\longrightarrow}g^{\ast}\mathcal{E}nd_{\mathcal{O}_Y}(\mathcal{E}). \]

Take the pull-back
\[
\xy
(-24,7.5)*+{V\times_YV=:}="8";
(-10,8)*+{V^{(2)}}="10";
(10,8)*+{V}="12";
(-10,-8)*+{V}="20";
(10,-8)*+{Y}="22";
(15,-9)*+{,}="23";
(0,0)*+{\square}="1";
{\ar^{q_2} "10";"12"};
{\ar_{q_1} "10";"20"};
{\ar_{g} "20";"22"};
{\ar^{g} "12";"22"};
\endxy
\]
and put
\begin{align*}
q^{(2)}&:=g\circ q_1=g\circ q_2,\\
\phi^{(2)}&:=(%
\mathcal{E}nd_{\mathcal{O}_{V^{(2)}}}(q^{(2)\ast}\mathcal{E})^{\times}%
\overset{\cong}{\rightarrow}%
q^{(2)\ast}\mathcal{E}nd_{\mathcal{O}_Y}(\mathcal{E})^{\times}%
\overset{q_1^{\ast}\phi^{-1}}{\longrightarrow}%
q^{(2)\ast}\mathcal{A}\\
&\qquad\overset{q_2^{\ast}\phi}{\longrightarrow}%
q^{(2)\ast}\mathcal{E}nd_{\mathcal{O}_Y}(\mathcal{E})^{\times}%
\overset{\cong}{\rightarrow}%
\mathcal{E}nd_{\mathcal{O}_{V^{(2)}}}(q^{(2)\ast}\mathcal{E})^{\times}%
).
\end{align*}
Then, since $q^{(2)\ast}\mathcal{E}nd_{\mathcal{O}_Y}(\mathcal{E})$ is an Azumaya algebra on $V^{(2)}$, there exists a surjective \'{e}tale morphism $W\rightarrow V^{(2)}$ and an element $c\in\Gamma(W,\mathcal{E}nd_{\mathcal{O}_Y}(q^{(2)\ast}\mathcal{E})\mid_W)^{\times}$
such that $\phi^{(2)}$ is the inner automorphism defined by $c$ :
\[ \phi^{(2)}\mid_W=\mathrm{Inn}(c) \]

By Assumption \ref{Assumption}, there exists a surjective \'{e}tale morphism $V^{\prime}\overset{g^{\prime}}{\longrightarrow}Y$ which factors through $V$
\[
\xy
(-12,-4)*+{V^{\prime}}="0";
(0,0)*+{V}="2";
(0,-9)*+{}="3";
(12,-4)*+{Y}="4";
{\ar^{v} "0";"2"};
{\ar^{g} "2";"4"};
{\ar@/_0.80pc/_{g^{\prime}} "0";"4"};
{\ar@{}|\circlearrowright "2";"3"};
\endxy
\]
such that the induced morphism
\[ V^{\prime(2)}:=V^{\prime}\times_YV^{\prime}\overset{v^{(2)}}{\longrightarrow}V^{(2)} \]
factors through $W$.

So, by replacing $V\overset{g}{\longrightarrow}Y$ by $V^{\prime}\overset{g^{\prime}}{\longrightarrow}Y$, we may assume the existence of a quartet
\[ (\mathcal{V},\mathcal{E},\phi,c) \]
which satisfies
\begin{align*}
\mathcal{V}&=(V\overset{g}{\longrightarrow}Y)
,\ \text{surjective \'{e}tale morphism of finite type},\\
\mathcal{E}&\in\mathrm{q\text{-}Coh}(Y),\ \text{locally free of finite rank},\\
\phi&:g^{\ast}\mathcal{A}\overset{\cong}{\longrightarrow}g^{\ast}\mathcal{E}nd_{\mathcal{O}_Y}(\mathcal{E}),\ \mathcal{O}_V\text{-algebra isomorphism},\\
c&\in\Gamma(V^{(2)},\mathcal{E}nd_{\mathcal{O}_Y}(q^{(2)\ast}\mathcal{E})^{\times}),\ \phi^{(2)}=\mathrm{Inn}(c).\\
\end{align*}
We call $(\mathcal{V},\mathcal{E},\phi,c)$ a {\it compatible trivialization} of $\mathcal{A}$.
Remark that for any refinement of $\mathcal{V}$
\[
\mathcal{V}^{\prime}=(V^{\prime}\overset{g^{\prime}}{\longrightarrow}V)
\qquad
\xy
(-12,0)*+{V^{\prime}}="0";
(0,4)*+{V}="2";
(0,-5)*+{}="3";
(12,0)*+{Y}="4";
(17,-1)*+{,}="5";
{\ar^{v} "0";"2"};
{\ar^{g} "2";"4"};
{\ar@/_0.80pc/_{g^{\prime}} "0";"4"};
{\ar@{}|\circlearrowright "2";"3"};
\endxy
\]
we obtain an induced compatible trivialization of $\mathcal{A}$ on $\mathcal{V}^{\prime}$
\[ (\mathcal{V}^{\prime},\mathcal{E},v^{\ast}\phi,v^{(2)\ast}c). \]

Let $q_{ij}:V^{(3)}=V\times_YV\times_YV\rightarrow V^{(2)}\ (1\le i<j\le 3)$ be the projections to the $(i,j)$-th components, and put $q^{(3)}:=q^{(2)}\circ q_{ij}$.
If we put
\[ \chi:=q_{12}^{\ast}c\cdot q_{13}^{\ast}c^{-1}\cdot q_{23}^{\ast}c\in\Gamma(V^{(3)},q^{(3)\ast}\mathcal{E}nd_{\mathcal{O}_Y}(\mathcal{E})^{\times}), \]
then in $\mathrm{Aut}_{\mathcal{O}_{V^{(3)}}}(q^{(3)\ast}\mathcal{E}nd_{\mathcal{O}_Y}(\mathcal{E}))$, we have
\begin{align*}
\mathrm{Inn}(\chi)&=\mathrm{Inn}(q_{12}^{\ast}c)\cdot\mathrm{Inn}(q_{13}^{\ast}c^{-1}\cdot)\mathrm{Inn}(q_{23}^{\ast}c)\\
&=q_{12}^{\ast}\phi^{(2)}\circ q_{13}^{\ast}\phi^{(2)-1}\circ q_{23}^{\ast}\phi^{(2)}\\
&=\mathrm{id}.
\end{align*}

So $\chi$ is in the center of $\Gamma(V^{(3)},q^{(3)\ast}\mathcal{E}nd_{\mathcal{O}_Y}(\mathcal{E})^{\times})$, i.e.,
\[ \chi\in Z(\Gamma(V^{(3)},q^{(3)\ast}\mathcal{E}nd_{\mathcal{O}_Y}(\mathcal{E})^{\times}))=\Gamma(V^{(3)},\mathcal{O}_{V^{(3)}}^{\times}). \]
Thus we obtain a 2-cocycle $\chi_Y^{\mathcal{V}}(\mathcal{A})$ in the \v{C}ech complex $\check{\mathcal{C}}^{\bullet}(\mathcal{V}/Y,\mathbb{G}_{m,Y})$, which defines $\chi_Y(\mathcal{A})\in \check{H}^2_{\mathrm{et}}(Y,\mathbb{G}_{m,Y})$.

For any Azumaya algebra $\mathcal{A}$ on $Y$, take a compatible trivialization $(\mathcal{V},\mathcal{E},\phi,c)$. By Lemma \ref{ec}, there exists $\mathcal{U}=(U\overset{f}{\longrightarrow}X)\in\mathrm{Cov}_{\mathrm{et}}(X)$ such that $\pi^{\ast}\mathcal{U}\le\mathcal{V}$. Here, $\pi^{\ast}\mathcal{U}\in\mathrm{Cov}_{\mathrm{et}}(Y)$ denotes the covering induced by pull-back by $\pi$.

So, replacing $(\mathcal{V},\mathcal{E},\phi,c)$ by the induced compatible trivialization on $\pi^{\ast}\mathcal{U}$, we may assume $\mathcal{V}=\pi^{\ast}\mathcal{U}$.

Take the pull-back
\[
\xy
(-30,8)*+{V^{(3)}}="8";
(-10,8)*+{V^{(2)}}="10";
(10,8)*+{V}="12";
(30,8)*+{Y}="14";
(-30,-8)*+{U^{(3)}}="18";
(-10,-8)*+{U^{(2)}}="20";
(8,-24)*+{}="21";
(10,-8)*+{U}="22";
(10,-14)*+{}="23";
(30,-8)*+{X}="24";
(-20,0)*+{\square}="1";
(0,0)*+{\square}="2";
(20,0)*+{\square}="3";
(52,4)*+{(\ell=1,2)}="4";
(56,-4)*+{(1\le i<j\le 3)}="5";
{\ar^{q_{\ell}} "10";"12"};
{\ar^{q_{ij}} "8";"10"};
{\ar_{\varpi^{(2)}} "10";"20"};
{\ar_{\varpi^{(3)}} "8";"18"};
{\ar^{p_{\ell}} "20";"22"};
{\ar^{p_{ij}} "18";"20"};
{\ar^{\varpi} "12";"22"};
{\ar^{g} "12";"14"};
{\ar^{f} "22";"24"};
{\ar^{\pi} "14";"24"};
{\ar@{}|\circlearrowright "20";"21"};
{\ar@{}|\circlearrowright "22";"23"};
{\ar@/_1.20pc/_{p^{(2)}} "20";"24"};
{\ar@/_2.40pc/_{p^{(3)}} "18";"24"};
\endxy
\]
Then we obtain a compatible trivialization $(\mathcal{U},\mathcal{N}_{\pi}\mathcal{E},N_{\varpi}\phi,N_{\varpi^{(2)}}(c))$ of $\mathcal{N}_{\pi}(\mathcal{A})$, defined as follows:
\begin{align*}
N_{\varpi}(\phi):=&(f^{\ast}\mathcal{N}_{\pi}\mathcal{A}\overset{\theta^f_{\pi}}{\longrightarrow}\mathcal{N}_{\varpi}g^{\ast}\mathcal{A}\overset{\mathcal{N}_{\varpi}(\phi)}{\longrightarrow}\mathcal{N}_{\varpi}g^{\ast}\mathcal{E}nd_{\mathcal{O}_Y}(\mathcal{E})\\
&\quad\overset{(\theta^f_{\pi})^{-1}}{\longrightarrow}f^{\ast}\mathcal{N}_{\pi}\mathcal{E}nd_{\mathcal{O}_Y}(\mathcal{E})\overset{f^{\ast}\delta_{\pi}}{\longrightarrow}f^{\ast}\mathcal{E}nd_{\mathcal{O}_X}(\mathcal{N}_{\pi}\mathcal{E})
\end{align*}
\[
\xy
(-36,6)*+{\Gamma(V^{(2)},\mathcal{E}nd_{\mathcal{O}_{V^{(2)}}}(q^{(2)\ast}\mathcal{E})^{\times})}="1";
(36,6)*+{\Gamma(U^{(2)},\mathcal{E}nd_{\mathcal{O}_{U^{(2)}}}(p^{(2)\ast}\mathcal{N}_{\pi}\mathcal{E})^{\times})}="2";
(-44,-6)*+{\mathrm{Aut}_{\mathcal{O}_{V^{(2)}}}(q^{(2)\ast}\mathcal{E})}="3";
(44,-6)*+{\mathrm{Aut}_{\mathcal{O}_{U^{(2)}}}(p^{(2)\ast}\mathcal{N}_{\pi}\mathcal{E})}="4";
(0,-20)*+{\mathrm{Aut}_{\mathcal{O}_{U^{(2)}}}(\mathcal{N}_{\varpi^{(2)}}q^{(2)\ast}\mathcal{E})}="5";
(0,12)*+{}="6";
{\ar^{N_{\varpi^{(2)}}} "1";"2"};
{\ar@{=} "1";"3"};
{\ar@{=} "2";"4"};
{\ar_{\mathcal{N}_{\varpi^{(2)}}} "3";"5"};
{\ar_{\cong}^<<<<<<<<{c(\theta^{p^{(2)}}_{\pi})} "4";"5"};
{\ar@{}|\circlearrowright "5";"6"};
\endxy
\]
Here $c(\theta^{p^{(2)}}_{\pi})$ is the conjugation by $\theta^{p^{(2)}}_{\pi}$.

By Remark \ref{IthetaRem}, there is an induced group isomorphism
\[ I_{\theta}^{\times}:\mathcal{E}nd_{\mathcal{O}_{U^{(2)}}}(\mathcal{N}_{\varpi^{(2)}}q^{(2)\ast}\mathcal{E})^{\times}\overset{\cong}{\longrightarrow}\mathcal{E}nd_{\mathcal{O}_{U^{(2)}}}(p^{(2)\ast}\mathcal{N}_{\pi}\mathcal{E})^{\times}. \]
\begin{claim}\label{FinClaim}
$(\mathcal{U},\mathcal{N}_{\pi}\mathcal{E},N_{\varpi}\phi,N_{\varpi^{(2)}}(c))$ is a compatible trivialization of $\mathcal{N}_{\pi}(\mathcal{A})$.
\end{claim}
\begin{pf}(Proof of Claim \ref{FinClaim})
It suffices to show
\begin{equation}
(N_{\varpi}(\phi))^{(2)}=\mathrm{Inn}(N_{\varpi^{(2)}}(c)).
\label{EqToShow}
\end{equation}

Using Proposition \ref{GeneralthetaProp} {\rm (ii)} and Proposition \ref{IthetaProp}, we can show easily
\begin{equation}
(N_{\varpi}(\phi))^{(2)}=I_{\theta}^{\times}\circ\delta_{\varpi^{(2)}}\circ\mathcal{N}_{\varpi^{(2)}}(\phi^{(2)})\circ\delta_{\varpi^{(2)}}^{-1}\circ I_{\theta}^{\times -1}.
\label{FinEq1}
\end{equation}
By Remark \ref{IthetaRem}, we have a commutative diagram 

\[
\xy
(-36,16)*+{\mathrm{Aut}_{\mathcal{O}_{U^{(2)}}}(\mathcal{N}_{\varpi^{(2)}}q^{(2)\ast}\mathcal{E})}="1";
(36,16)*+{\mathrm{Aut}_{\mathcal{O}_{U^{(2)}}}(p^{(2)\ast}\mathcal{N}_{\pi}\mathcal{E})}="2";
(-36,0)*+{\Gamma(U^{(2)},\mathcal{E}nd_{\mathcal{O}_{U^{(2)}}}(\mathcal{N}_{\varpi^{(2)}}q^{(2)\ast}\mathcal{E})^{\times})}="11";
(36,0)*+{\Gamma(U^{(2)},\mathcal{E}nd_{\mathcal{O}_{U^{(2)}}}(p^{(2)\ast}\mathcal{N}_{\pi}\mathcal{E})^{\times})}="12";
(-36,-16)*+{\mathrm{Aut}_{\mathcal{O}_{U^{(2)}}}(\mathcal{E}nd_{\mathcal{O}_{U^{(2)}}}(\mathcal{N}_{\varpi^{(2)}}q^{(2)\ast}\mathcal{E})^{\times})}="21";
(36,-16)*+{\mathrm{Aut}_{\mathcal{O}_{U^{(2)}}}(\mathcal{E}nd_{\mathcal{O}_{U^{(2)}}}(p^{(2)\ast}\mathcal{N}_{\pi}\mathcal{E})^{\times})\, .}="22";
{\ar^{c(\theta^{p^{(2)}}_{\pi})^{-1}} "1";"2"};
{\ar^{\Gamma(U^{(2)},I_{\theta}^{\times})} "11";"12"};
{\ar@{=} "1";"11"};
{\ar@{=} "2";"12"};
{\ar^{\mathrm{Inn}} "11";"21"};
{\ar^{\mathrm{Inn}} "12";"22"};
{\ar_{c(I_{\theta}^{\times})} "21";"22"};
{\ar@{}|\circlearrowright "1";"12"};
{\ar@{}|\circlearrowright "11";"22"};
\endxy
\]

Thus we obtain the following commutative diagram:
\[
\xy
(-32,12)*+{\mathrm{Aut}_{\mathcal{O}_{V^{(2)}}}(q^{(2)\ast}\mathcal{E})}="0";
(-32,0)*+{\mathrm{Aut}_{\mathcal{O}_{U^{(2)}}}(\mathcal{N}_{\varpi^{(2)}}q^{(2)\ast}\mathcal{E})}="1";
(-32,-12)*+{\mathrm{Aut}_{\mathcal{O}_{U^{(2)}}}(p^{(2)\ast}\mathcal{N}_{\pi}\mathcal{E})}="2";
(32,24)*+{\mathrm{Aut}_{\mathcal{O}_{V^{(2)}}}(\mathcal{E}nd_{\mathcal{O}_{V^{(2)}}}(q^{(2)\ast}\mathcal{E})^{\times})}="3";
(32,8)*+{\mathrm{Aut}_{\mathcal{O}_{U^{(2)}}}(\mathcal{N}_{\varpi^{(2)}}\mathcal{E}nd_{\mathcal{O}_{V^{(2)}}}(q^{(2)\ast}\mathcal{E})^{\times})}="4";
(32,-8)*+{\mathrm{Aut}_{\mathcal{O}_{U^{(2)}}}(\mathcal{E}nd_{\mathcal{O}_{U^{(2)}}}(\mathcal{N}_{\varpi^{(2)}}q^{(2)\ast}\mathcal{E})^{\times})}="5";
(32,-24)*+{\mathrm{Aut}_{\mathcal{O}_{U^{(2)}}}(\mathcal{E}nd_{\mathcal{O}_{U^{(2)}}}(p^{(2)\ast}\mathcal{N}_{\pi}\mathcal{E})^{\times})}="6";
(-32,24)*+{\Gamma(V^{(2)},\mathcal{E}nd_{\mathcal{O}_{V^{(2)}}}(q^{(2)\ast}\mathcal{E})^{\times})}="7";
(-32,-24)*+{\Gamma(U^{(2)},\mathcal{E}nd_{\mathcal{O}_{U^{(2)}}}(p^{(2)\ast}\mathcal{N}_{\pi}\mathcal{E})^{\times})}="8";
{\ar@{=} "0";"7"};
{\ar@{=} "2";"8"};
{\ar_{\mathcal{N}_{\varpi^{(2)}}} "0";"1"};
{\ar_{c(\theta^{p^{(2)}}_{\pi})^{-1}} "1";"2"};
{\ar^{\mathrm{Inn}} "7";"3"};
{\ar_{\mathrm{Inn}} "8";"6"};
{\ar_{\mathcal{N}_{\varpi^{(2)}}} "3";"4"};
{\ar_{c(\delta_{\varpi^{(2)}})} "4";"5"};
{\ar_{c(I_{\theta}^{\times})} "5";"6"};
{\ar@{}|\circlearrowright "7";"6"};
\endxy
\]

By $(\ref{FinEq1})$, this means $(\ref{EqToShow})$.
\end{pf}

We have
\[ \check{H}^2_{\mathrm{et}}(\mathcal{V}/Y,\mathbb{G}_{m,Y})=\check{H}^2_{\mathrm{et}}(\pi^{\ast}\mathcal{U}/Y,\mathbb{G}_{m,Y})=\check{H}^2_{\mathrm{et}}(\mathcal{U}/X,\pi_{\ast}\mathbb{G}_{m,Y}), \]
and the canonical natural isomorphism
\[ \mathfrak{c}^{-1}:H^2_{\mathrm{et}}(Y,\mathbb{G}_{m,Y})\overset{\cong}{\longrightarrow}H^2_{\mathrm{et}}(X,\pi_{\ast}\mathbb{G}_{m,Y}) \]
fits into the following commutative diagram:
\[
\xy
(-22,8)*+{\check{H}^2_{\mathrm{et}}(\mathcal{V}/Y,\mathbb{G}_{m,Y})}="0";
(22,8)*+{\check{H}^2_{\mathrm{et}}(\mathcal{U}/X,\pi_{\ast}\mathbb{G}_{m,Y})}="2";
(-22,-8)*+{H^2_{\mathrm{et}}(Y,\mathbb{G}_{m,Y})}="4";
(22,-8)*+{H^2_{\mathrm{et}}(X,\pi_{\ast}\mathbb{G}_{m,X})}="6";
{\ar^{\mathrm{id}}_{=} "0";"2"};
{\ar^{\cong}_{\mathfrak{c}^{-1}} "4";"6"};
{\ar_{\mathrm{can.}} "0";"4"};
{\ar^{\mathrm{can.}} "2";"6"};
{\ar@{}|\circlearrowright "0";"6"};
\endxy
\]
So, it suffices to show
\[ \check{H}^2_{\mathrm{et}}(\mathcal{U}/X,N_{\pi})(\chi^{\mathcal{V}}_Y(\mathcal{A}))=\chi^{\mathcal{U}}_X(\mathcal{N}_{\pi}(\mathcal{A})). \]

Similarly as $N_{\varpi^{(2)}}$, we can construct a homomorphism
\[ N_{\varpi^{(3)}}:\Gamma(V^{(3)},q^{(3)\ast}\mathcal{E}nd_{\mathcal{O}_Y}(\mathcal{E})^{\times})\rightarrow\Gamma(U^{(3)},p^{(3)\ast}\mathcal{E}nd_{\mathcal{O}_X}(\mathcal{N}_{\pi}\mathcal{E})^{\times}), \]
compatible with $N_{\varpi^{(2)}}$ and $N_{\varpi^{(3)}}:\Gamma(V^{(3)},\mathcal{O}_{V^{(3)}}^{\times})\rightarrow\Gamma(U^{(3)},\mathcal{O}_{U^{(3)}}^{\times})$.
\[
\xy
(-32,24)*+{\Gamma(V^{(2)},\mathcal{E}nd_{\mathcal{O}_{V^{(2)}}}(q^{(2)\ast}\mathcal{E})^{\times})}="10";
(-32,12)*+{\Gamma(V^{(2)},q^{(2)\ast}\mathcal{E}nd_{\mathcal{O}_Y}(\mathcal{E})^{\times})}="12";
(-32,0)*+{\Gamma(V^{(3)},q^{(3)\ast}\mathcal{E}nd_{\mathcal{O}_Y}(\mathcal{E})^{\times})}="14";
(-32,-12)*+{Z(\Gamma(V^{(3)},q^{(3)\ast}\mathcal{E}nd_{\mathcal{O}_Y}(\mathcal{E})^{\times}))}="16";
(-32,-24)*+{\Gamma(V^{(3)},\mathcal{O}_{V^{(3)}}^{\times})}="18";
(32,24)*+{\Gamma(U^{(2)},\mathcal{E}nd_{\mathcal{O}_{U^{(2)}}}(p^{(2)\ast}\mathcal{N}_{\pi}\mathcal{E})^{\times})}="20";
(32,12)*+{\Gamma(U^{(2)},p^{(2)\ast}\mathcal{E}nd_{\mathcal{O}_X}(\mathcal{N}_{\pi}\mathcal{E})^{\times})}="22";
(32,0)*+{\Gamma(U^{(3)},p^{(3)\ast}\mathcal{E}nd_{\mathcal{O}_X}(\mathcal{N}_{\pi}\mathcal{E})^{\times})}="24";
(32,-12)*+{Z(\Gamma(U^{(3)},p^{(3)\ast}\mathcal{E}nd_{\mathcal{O}_X}(\mathcal{N}_{\pi}\mathcal{E})^{\times}))}="26";
(32,-24)*+{\Gamma(U^{(3)},\mathcal{O}_{U^{(3)}}^{\times})}="28";
{\ar^{N_{\varpi^{(2)}}} "10";"20"};
{\ar^{\exists N_{\varpi^{(3)}}} "14";"24"};
{\ar_{N_{\varpi^{(3)}}} "18";"28"};
{\ar^{\cong} "10";"12"};
{\ar^{\cong} "20";"22"};
{\ar^{q_{ij}^{\ast}} "12";"14"};
{\ar^{p_{ij}^{\ast}} "22";"24"};
{\ar@{_(->} "16";"14"};
{\ar@{_(->} "26";"24"};
{\ar^{\cong} "18";"16"};
{\ar^{\cong} "28";"26"};
{\ar@{}|\circlearrowright "10";"24"};
{\ar@{}|\circlearrowright "14";"28"};
\endxy
\]

From this, we have
\begin{align*}
\check{H}^2_{\mathrm{et}}(\mathcal{U}/X,N_{\pi})(\chi^{\mathcal{V}}_Y(\mathcal{A}))&=N_{\varpi^{(3)}}(q_{12}^{\ast}c\cdot q_{13}^{\ast}c^{-1}\cdot q_{23}^{\ast}c)\\
&=N_{\varpi^{(3)}}(q_{12}^{\ast}c)\cdot N_{\varpi^{(3)}}(q_{13}^{\ast}c^{-1})\cdot N_{\varpi^{(3)}}(q_{23}^{\ast}c)\\
&=(p_{12}^{\ast}N_{\varpi^{(2)}}(c))\cdot (p_{13}^{\ast}N_{\varpi^{(2)}}(c)^{-1})\cdot (p_{23}^{\ast}N_{\varpi^{(2)}}(c))\\
&=\chi^{\mathcal{U}}_X(\mathcal{N}_{\pi}(\mathcal{A})).
\end{align*}
\end{pf}

\section{Brauer-Mackey functor on the Galois category}

Let $\mathrm{Ab}$ be the category of abelian groups.
For any profinite group $G$, let $G\mathrm{\text{-}Sp}$ denote the category of finite discrete $G$-spaces and continuous equivariant $G$-maps.
\begin{defn}
\label{DefOfMackFtr}
Let $\mathcal{C}$ be a 
Galois category, with fundamental functor $F$. In other words, there exists a profinite group $\pi(\mathcal{C})$ such that $F$ gives an equivalence from $\mathcal{C}$ to $\pi(\mathcal{C})\text{-}\mathrm{Sp}$. $($For the precise definition of Galois category, see {\rm \cite{Murre}}$)$.

A cohomological Mackey functor on $\mathcal{C}$ is a pair of functors $M=(M^{\ast},M_{\ast})$ from $\mathcal{C}$ to $\mathrm{Ab}$, where $M^{\ast}$ is contravariant and $M_{\ast}$ is covariant,
satisfying the following conditions $:$\\
{\rm (0)} $M^{\ast}(X)=M_{\ast}(X)(=:M(X))\quad(\forall X\in\mathrm{Ob}(\mathcal{C}))$.\\
{\rm (1)} {\rm (Additivity)} For each coproduct $X\overset{i_X}{\hookrightarrow}X\coprod Y\overset{i_Y}{\hookleftarrow}Y$ in $\mathcal{C}$, canonical morphism
\[ (M^{\ast}(i_X),M^{\ast}(i_Y)):M(X\textstyle{\coprod}Y)
{\rightarrow}M(X)\oplus M(Y) \]
is an isomorphism.\\
{\rm (2)} {\rm (Mackey condition)} For any pull-back diagram
\[
\xy
(-8,6)*+{Y^{\prime}}="0";
(-8,-6)*+{X^{\prime}}="2";
(8,6)*+{Y}="4";
(8,-6)*+{X}="6";
(12,-7)*+{,}="7";
(0,0)*+{\square}="8";
{\ar_{\pi} "0";"2"};
{\ar^{\varpi^{\prime}} "0";"4"};
{\ar^{\pi^{\prime}} "4";"6"};
{\ar_{\varpi} "2";"6"};
\endxy
\]
the following diagram is commutative $:$
\[
\xy
(-12,7)*+{M(Y)}="0";
(-12,-7)*+{M(X)}="2";
(12,7)*+{M(Y^{\prime})}="4";
(12,-7)*+{M(X^{\prime})}="6";
{\ar_{M_{\ast}(\pi)} "0";"2"};
{\ar^{M^{\ast}(\varpi^{\prime})} "0";"4"};
{\ar^{M_{\ast}(\pi^{\prime})} "4";"6"};
{\ar_{M^{\ast}(\varpi)} "2";"6"};
{\ar@{}|\circlearrowright "0";"6"};
\endxy
\]\\
{\rm (3)} {\rm (Cohomological condition)}
For any morphism $\pi:X\rightarrow Y$ in $\mathcal{C}$ with $X$ and $Y$ connected $($ i.e. not decomposable into non-trivial coproducts $)$, we have
\[ M_{\ast}(\pi)\circ M^{\ast}(\pi)=\text{multiplication by}\ \deg\pi \]
where $\deg\pi:=\sharp F(Y)/\sharp F(X)$.
\[ 
\xy
(-24,0)*+{M(X)}="0";
(0,8)*+{M(Y)}="2";
(24,0)*+{M(X)}="4";
(0,-4)*+{}="6";
{\ar^{M^{\ast}(\pi)} "0";"2"};
{\ar^{M_{\ast}(\pi)} "2";"4"};
{\ar@/_0.60pc/_{\deg \pi} "0";"4"};
{\ar@{}|\circlearrowright "2";"6"};
\endxy
\]
\end{defn}

\begin{defn}
Let $M$ and $N$ be Mackey functors on $\mathcal{C}$.
A morphism $f:M\rightarrow N$ is a collection of homomorphisms in $\mathrm{Ab}$
\[ \{ f(X)\mid X\in\mathcal{C} \} , \]
which is natural with respect to each of the covariant and the contravariant part of $M$ and $N$.
With the objectwise composition, we define the category of cohomological Mackey functors $\mathrm{Mack}_c(\mathcal{C})$.
\end{defn}

A standard example is the cohomological Mackey functor on a profinite group (see \cite{Bley-Boltje}):
\begin{defn}
Let $G$ be a profinite group, and let $\mathcal{C}=G\text{-}\mathrm{Sp}$, $F=\mathrm{id}_{\mathcal{C}}$. A cohomological Mackey functor on $\mathcal{C}$ is simply called a cohomological Mackey functor on $G$, and their category is denoted by $Mack_c(G)$.
\end{defn}
\begin{rem}\label{Rem6.3}
Any object $X$ in $G\text{-}\mathrm{Sp}$ is a direct sum of transitive $G$-sets of the form $G/H$, where $H$ is a open subgroup of $G$. So a Mackey functor on $G$ is equal to the following datum $:$\\
- an abelian group $M(H)$ for each open $H\le G$, with structure maps:\\
\ \ - a homomorphism $\mathrm{res}^H_K:M(H)\rightarrow M(K)$ 
for each open $K\le H\le G$,\\
\ \ - a homomorphism $\mathrm{cor}^H_K:M(K)\rightarrow M(H)$ 
for each open $K\le H\le G$,\\
\ \ - a homomorphism $c_{g,H}:M(H)\rightarrow M({}^g\! H)$ 
for each open $H\le G$ and $g\in G$,\\
where ${}^g\! H:=gHg^{-1}$, satisfying certain compatibilities $($cf. {\rm \cite{Bley-Boltje}}$)$.
Here, $M(G/H)$ is abbreviated to $M(H)$ for any open subgroup $H\le G$.
\end{rem}

\begin{defn}
Let $G$ be a finite group, and let $G^{\mathrm{op}}$ be its opposite group. For any Mackey functor $M=(M,\mathrm{res},\mathrm{cor},c)\in\mathrm{Mack}_c(G)$ $($in the notation of Remark \ref{Rem6.3}$)$, we define its {\it opposite Mackey functor} $M^{\mathrm{op}}$ by
\begin{align*}
M^{\mathrm{op}}(H^{\mathrm{op}})&:=M(H)\qquad(H\le G)\\
\mathrm{res}^{H^{\mathrm{op}}}_{K^{\mathrm{op}}}&:=\mathrm{res}^H_K\ \ \,\qquad(K\le H\le G)\\
\mathrm{cor}^{H^{\mathrm{op}}}_{K^{\mathrm{op}}}&:=\mathrm{cor}^H_K\ \ \qquad(K\le H\le G)\\
c_{g,H^{\mathrm{op}}}&:=c_{g^{-1},H}\ \qquad (g\in G,H\le G).
\end{align*}
This gives an isomorphism of categories
\[ \mathrm{op}:\mathrm{Mack}_c(G)\rightarrow\mathrm{Mack}_c(G^{\mathrm{op}}). \]
\end{defn}


For any finite \'{e}tale covering $\pi:Y\rightarrow X$, put $\mathrm{Br}^{\ast}(\pi):=\pi^{\ast}$ and $\mathrm{Br}_{\ast}(\pi):=N_{\pi}$.
Then we obtain a cohomological Mackey functor $\mathrm{Br}$ (and similarly $\mathrm{Br}^{\prime}$, $H^2_{\mathrm{et}}(-,\mathbb{G}_m)$) as follows. Remark that for any connected scheme $S$, the category $(\mathrm{FEt}/S)$ of finite \'{e}tale coverings over $S$ becomes a Galois category \cite{Murre}.
\begin{thm}\label{MainThm}
For any connected scheme $S$ satisfying Assumption \ref{Assumption}, we have a sequence of cohomological Mackey functors on $(\mathrm{FEt}/S)$
\[ \mathrm{Br}\hookrightarrow \mathrm{Br}^{\prime}\hookrightarrow H^2_{\mathrm{et}}(-,\mathbb{G}_m).\]
\end{thm}
\begin{pf}
We only show Mackey and cohomological conditions. Since $\pi^{\ast}$ and $N_{\pi}$ are compatible with inclusions
\[ \mathrm{Br}(X)\hookrightarrow \mathrm{Br}^{\prime}(X)\hookrightarrow H^2_{\mathrm{et}}(X,\mathbb{G}_{m,X}),\]
it suffices to show for $H^2_{\mathrm{et}}(-, \mathbb{G}_m)$.

\noindent\underline{Mackey condition}

Let
\begin{equation}
\xy
(-8,6)*+{Y}="0";
(-8,-6)*+{X}="2";
(8,6)*+{Y^{\prime}}="4";
(8,-6)*+{X^{\prime}}="6";
(0,0)*+{\square}="8";
{\ar_{\pi} "0";"2"};
{\ar_{\varpi^{\prime}} "4";"0"};
{\ar^{\pi^{\prime}} "4";"6"};
{\ar^{\varpi} "6";"2"};
\endxy
\label{MackCondPrf}
\end{equation}
be a pull-back diagram in $(\mathrm{FEt}/S)$.

For any \'{e}tale morphism of finite type $f:U\rightarrow X$, take the pull-back of $(\ref{MackCondPrf})$ by $f$:
\[
\xy
(-8,6)*+{V}="0";
(-8,-6)*+{U}="2";
(8,6)*+{V^{\prime}}="4";
(8,-6)*+{U^{\prime}}="6";
(0,0)*+{\square}="8";
{\ar_{\pi_U} "0";"2"};
{\ar_{\varpi^{\prime}_U} "4";"0"};
{\ar^{\pi^{\prime}_U} "4";"6"};
{\ar^{\varpi_U} "6";"2"};
\endxy
\]
Then we have a commutative diagram
\[
\xy
(-32,0)*+{\mathrm{Aut}_{\mathcal{O}_V}(\mathcal{O}_V)}="0";
(-24,12)*+{\mathrm{Aut}_{\mathcal{O}_{V^{\prime}}}(\varpi_U^{\prime\ast}\mathcal{O}_V)}="2";
(16,8)*+{\mathrm{Aut}_{\mathcal{O}_{U^{\prime}}}(\mathcal{N}_{\pi_U^{\prime}}\varpi_U^{\prime\ast}\mathcal{O}_V)}="4";
(-24,-12)*+{\mathrm{Aut}_{\mathcal{O}_U}(\mathcal{N}_{\pi_U}\mathcal{O}_V)}="6";
(20,0)*+{}="5";
(16,-8)*+{\mathrm{Aut}_{\mathcal{O}_{U^{\prime}}}(\varpi_U^{\ast}\mathcal{N}_{\pi_U}\mathcal{O}_V)}="8";
(52,0)*+{\mathrm{Aut}_{\mathcal{O}_{U^{\prime}}}(\mathcal{O}_{U^{\prime}})}="10";
(8,0)*+{}="11";
{\ar_{\mathcal{N}_{\pi_U}} "0";"6"};
{\ar^{\mathcal{N}_{\pi_U^{\prime}}} "2";"4"};
{\ar^{\varpi_U^{\prime\ast}} "0";"2"};
{\ar_{\varpi_U^{\ast}} "6";"8"};
{\ar_{c(\theta^{\varpi_U}_{\pi_U})} "4";"8"};
{\ar^{\cong} "4";"10"};
{\ar_{\cong} "8";"10"};
{\ar@{}|\circlearrowright "0";"5"};
{\ar@{}|\circlearrowright "10";"11"};
\endxy
\]
where $c(\theta^{\varpi_U}_{\pi_U})$ is the conjugation by $\theta^{\varpi_U}_{\pi_U}$.

Thus we have a commutative diagram
\[
\xy
(-18,10)*+{\pi_{\ast}\mathbb{G}_{m,Y}}="0";
(-18,-8)*+{\mathbb{G}_{m,X}}="2";
(13,10)*+{\pi_{\ast}\varpi^{\prime}_{\ast}\mathbb{G}_{m,Y^{\prime}}}="4";
(34,5)*+{\varpi_{\ast}\pi^{\prime}_{\ast}\mathbb{G}_{m,Y^{\prime}}}="5";
(34,-8)*+{\varpi_{\ast}\mathbb{G}_{m,X^{\prime}}}="6";
(23.5,7.5)*+{\rotatebox{-24}{$\cong$}}="8";
(43,-9)*+{.}="10";
{\ar_{N_{\pi}} "0";"2"};
{\ar^{\pi_{\ast}(\varpi^{\prime}_{\sharp})} "0";"4"};
{\ar^{\varpi_{\ast}(N_{\pi^{\prime}})} "5";"6"};
{\ar_{\varpi_{\sharp}} "2";"6"};
{\ar@{}|\circlearrowright "0";"6"};
\endxy
\]

This yields a commutative diagram
\[
\xy
(-18,8)*+{H^2_{\mathrm{et}}(Y,\mathbb{G}_{m,Y})}="0";
(-18,-8)*+{H^2_{\mathrm{et}}(X,\mathbb{G}_{m,X})}="2";
(18,8)*+{H^2_{\mathrm{et}}(Y^{\prime},\mathbb{G}_{m,Y^{\prime}})}="4";
(18,-8)*+{H^2_{\mathrm{et}}(X^{\prime},\mathbb{G}_{m,X^{\prime}})}="6";
(32,-9)*+{.}="10";
{\ar_{N_{\pi}} "0";"2"};
{\ar^{\varpi^{\prime\ast}} "0";"4"};
{\ar^{N_{\pi^{\prime}}} "4";"6"};
{\ar_{\varpi^{\ast}} "2";"6"};
{\ar@{}|\circlearrowright "0";"6"};
\endxy
\]

\noindent\underline{Cohomological condition}

For any finite \'{e}tale covering $\varpi:V\rightarrow U$ of constant degree $d$, the composition
\begin{align*} \mathrm{Aut}_{\mathcal{O}_U}(\mathcal{O}_U)\overset{\varpi^{\ast}}{\rightarrow}\mathrm{Aut}_{\mathcal{O}_V}(\varpi^{\ast}\mathcal{O}_U)&\overset{\cong}{\rightarrow}\mathrm{Aut}_{\mathcal{O}_V}(\mathcal{O}_V)\\
&\overset{\mathcal{N}_{\varpi}}{\rightarrow}\mathrm{Aut}_{\mathcal{O}_U}(\mathcal{N}_{\varpi}\mathcal{O}_V)\overset{\cong}{\rightarrow}\mathrm{Aut}_{\mathcal{O}_U}(\mathcal{O}_U)
\end{align*}
is equal to the multiplication by $d$. This follows from the trivial case $\nabla:\underset{d}{\coprod}U\rightarrow U$ via fpqc descent.

From this, we can see
\[ N_{\pi}\circ \pi_{\sharp}:\mathbb{G}_{m,X}\rightarrow\mathbb{G}_{m,X} \]
is equal to the multiplication by $d=\deg (\pi)$
\[
\xy
(-24,0)*+{\mathbb{G}_{m,X}}="0";
(0,4)*+{\pi_{\ast}\mathbb{G}_{m,Y}}="2";
(24,0)*+{\mathbb{G}_{m,X}}="4";
(29,-2)*+{.}="5";
(0,-4)*+{}="6";
{\ar^{\pi_{\sharp}} "0";"2"};
{\ar^{N_{\pi}} "2";"4"};
{\ar@/_0.60pc/_{d} "0";"4"};
{\ar@{}|\circlearrowright "2";"6"};
\endxy
\]
Thus we obtain $N_{\pi}\circ\pi^{\ast}=d$.
\[
\xy
(-42,-4)*+{H^2_{\mathrm{et}}(X,\mathbb{G}_{m,X})}="0";
(0,10)*+{H^2_{\mathrm{et}}(Y,\mathbb{G}_{m,Y})}="2";
(42,-4)*+{H^2_{\mathrm{et}}(X,\mathbb{G}_{m,X})}="4";
(0,-2)*+{H^2_{\mathrm{et}}(X,\pi_{\ast}\mathbb{G}_{m,Y})}="6";
(12,6)*+{}="8";
(-12,6)*+{}="10";
(0,-17)*+{}="12";
{\ar^{\pi^{\ast}} "0";"2"};
{\ar^{N_{\pi}} "2";"4"};
{\ar^{\cong}_{\mathfrak{c}} "6";"2"};
{\ar@{}|\circlearrowright "0";"8"};
{\ar@{}|\circlearrowright "4";"10"};
{\ar_{H^2_{\mathrm{et}}(\pi_{\sharp})} "0";"6"};
{\ar_{H^2_{\mathrm{et}}(N_{\pi})} "6";"4"};
{\ar@/_2.20pc/_{d} "0";"4"};
{\ar@{}|\circlearrowright "6";"12"};
\endxy
\]
\end{pf}

\section{Restriction to a finite Galois covering}
Thus we have obtained a cohomological Mackey functor $\mathrm{Br}$ on $\mathrm{FEt}/S$. Pulling back by a quasi-inverse $\mathcal{S}$ of the fundamental functor
\[F:\mathrm{FEt}/S\overset{\simeq}{\longrightarrow}\pi(S)\text{-}\mathrm{Sp},\]
we obtain a Mackey functor on  $\pi(S)\,:$
\begin{cor}\label{Cor7.1}
There is a sequence of cohomological Mackey functors
\[ \mathrm{Br}\circ\mathcal{S}\hookrightarrow \mathrm{Br}^{\prime}\circ\mathcal{S}\hookrightarrow H^2_{\mathrm{et}}(-,\mathbb{G}_m)\circ\mathcal{S} \]
on $\pi(S)$, where $\mathrm{Br}\circ\mathcal{S}:=(\mathrm{Br}^{\ast}\circ\mathcal{S},\mathrm{Br}_{\ast}\circ\mathcal{S})$ $($and similarly for $\mathrm{Br}^{\prime}$, $H^2_{\mathrm{et}}(-,\mathbb{G}_m)\,)$.
\end{cor}

\begin{cor}\label{MainCor1}
Let $X$ be a connected scheme satisfying Assumption \ref{Assumption}. For any finite Galois covering $\pi:Y\rightarrow X$ with $\mathrm{Gal}(Y/X)=G$, there exists a cohomological Mackey functor $\BBr$ on $G$ which satisfies
\[ \BBr(H)\cong \mathrm{Br}(Y/H)\quad (\forall H\le G),\]
with structure maps $\varpi^{\ast}$ and $N_{\varpi}$ for each intermediate covering $\varpi$. $($We abbreviate $\BBr(G/H)$ to $\BBr(H)$, as in Remark \ref{Rem6.3}.$)$
\end{cor}
\begin{pf}
By the projection $\mathrm{pr}:\pi(X)\rightarrow\!\!\!\!\rightarrow G^{\mathrm{op}}$, we can regard any finite $G^{\mathrm{op}}$-set naturally as a finite $\pi(X)$-space, to obtain a functor
\[ G^{\mathrm{op}}\text{-}\mathrm{Sp}\rightarrow \pi(X)\text{-}\mathrm{Sp}. \]
Pulling back by this functor, and taking the opposite Mackey functor,\\
we obtain
\[
\xy
(-28,4)*+{\mathrm{Mack}_c(\pi(X))}="-2";
(2,4)*+{\mathrm{Mack}_c(G^{\mathrm{op}})}="0";
(30,4)*+{\mathrm{Mack}_c(G)}="2";
(-28,0)*+{\rotatebox{90}{$\in$}  }="4";
(30,0)*+{\rotatebox{90}{$\in$}  }="6";
(-28,-4)*+{M}="8";
(30,-4)*+{M_G}="10";
(34,-5)*+{.}="12";
{\ar "-2";"0"};
{\ar^{\mathrm{op}}"0";"2"};
{\ar@{|->} "8";"10"};
\endxy
\]

In terms of subgroups of $G$, $M_G$ satisfies
\[ M_G(H)=M(\mathrm{pr}^{-1}(H^{\mathrm{op}}))\quad(\forall H\le G).\]

Applying this to $\mathrm{Br}\circ\mathcal{S}$, we obtain $\BBr:=(Br\circ\mathcal{S})_G\ \in \mathrm{Mack}_c(G)$. Since the equivalence $\mathcal{S}:
(\pi(X)\text{-}\mathrm{Sp})\overset{\simeq}{\longrightarrow}(\mathrm{FEt}/X)$ satisfies
\[ \mathcal{S}(\pi(X)/\mathrm{pr}^{-1}(H^{\mathrm{op}}))\cong Y/H, \]
we have
\[ \BBr(H)\cong \mathrm{Br}(Y/H). \]
\end{pf}

\begin{cor}
Let $\pi:Y\rightarrow X$ be a finite Galois covering of a connected scheme $X$ satisfying Assumption \ref{Assumption}, with Galois group $G$. By a similar way, we can define $\BBr^{\prime}$ $($and also $(H^2_{\mathrm{et}}(-,\mathbb{G}_m)\circ\mathcal{S})_G\,)$.

Since $\mathrm{Mack}_c(G)$ is an abelian category with objectwise $($co-$)$kernels $($see for example {\rm \cite{Bouc}}$)$, we can take the quotient Mackey functor $\BBr^{\prime}/\BBr\in\mathrm{Mack}_c(G)$, which satisfies
\[ (\BBr^{\prime}/\BBr)(H)\cong(\mathrm{Br}^{\prime}(Y/H))/(\mathrm{Br}(Y/H)). \]
\end{cor}

\section{Appendix}

\subsection{Application of Bley and Boltje's theorem}

Let $\ell$ be a prime number. For any abelian group $A$, let
\[  A(\ell):=\{ m\in A\mid \exists e\in \mathbb{N}_{\ge 0},\ell^{e}m=0\} \]
be the $\ell$-primary part. This is a $\mathbb{Z}_{\ell}$-module.

\begin{defn}[\cite{Bley-Boltje}]
For any finite group $H$,\\
$H$ is {\it $\ell$-hypoelementary} $\underset{\text{def}}{\Leftrightarrow}$ $H$ has a normal $\ell$-subgroup with a cyclic quotient.\\
$H$ is {\it hypoelementary} $\underset{\text{def}}{\Leftrightarrow}$ $H$ is $\ell$-hypoelementary for some prime $\ell$.
\end{defn}

\begin{fact}[\cite{Bley-Boltje}]
\label{Bley-Boltje's Thm}
Let $M$ be a cohomological Mackey functor on a finite group $G$.\\
{\rm (i)} Let $\ell$ be a prime number. If $H\le G$ is not $\ell$-{\it hypoelementary}, then there is a natural isomorphism of $\mathbb{Z}_{\ell}$-modules
\[ \bigoplus _{\underset{n:\text{odd}}{U=H_0<\cdots<H_n=H}}M(U)(\ell)^{|U|}\cong\bigoplus_{\underset{n:\text{even}}{U=H_0<\cdots<H_n=H}}M(U)(\ell)^{|U|}. \]\\
{\rm (ii)} If $H\le G$ is not {\it hypoelementary} and $M(U)$ is torsion for any subgroup $U\le H$, then there is a natural isomorphism of abelian groups
\[ \bigoplus _{\underset{n:\text{odd}}{U=H_0<\cdots<H_n=H}}M(U)^{|U|}\cong\bigoplus_{\underset{n:\text{even}}{U=H_0<\cdots<H_n=H}}M(U)^{|U|}. \]
Here, $|U|$ denotes the order of $U$.
\end{fact}

Applying this theorem to $\BBr$, we obtain the following relations for the Brauer groups of intermediate coverings:

\begin{cor}
\label{MainCor2}
Let $X$ be a connected scheme satisfying Assumption \ref{Assumption} and $\pi:Y\rightarrow X$ be a finite Galois covering with $\mathrm{Gal}(Y/X)=G$.\\
{\rm (i)} Let $\ell$ be a prime number. If $H\le G$ is not $\ell$-hypoelementary, then there is a natural isomorphism of $\mathbb{Z}_{\ell}$-modules
\[ \bigoplus _{\underset{n:\text{odd}}{U=H_0<\cdots<H_n=H}}\mathrm{Br}(Y/U)(\ell)^{|U|}\cong\bigoplus_{\underset{n:\text{even}}{U=H_0<\cdots<H_n=H}}\mathrm{Br}(Y/U)(\ell)^{|U|}. \]\\
{\rm (ii)} If $H\le G$ is not hypoelementary, then there is a natural isomorphism of abelian groups
\[ \bigoplus _{\underset{n:\text{odd}}{U=H_0<\cdots<H_n=H}}\mathrm{Br}(Y/U)^{|U|}\cong\bigoplus_{\underset{n:\text{even}}{U=H_0<\cdots<H_n=H}}\mathrm{Br}(Y/U)^{|U|}. \]
\end{cor}

Finally, we derive some numerical equations related to Brauer groups from Corollary \ref{MainCor2}.

\begin{defn}
Let $G$ be a finite group. For any subgroups $U\le H\le G$, put
\[\mu(U,H):=\sum_{U=H_0<\cdots <H_n=H}(-1)^n,\quad \text{M\"{o}bius function}.\]\end{defn}

If $m$ (resp. $m_{\ell}$) is an additive invariant of abelian groups (resp. $\mathbb{Z}_{\ell}$-modules) which is finite on Brauer groups, we obtain the following equations:
\begin{cor}
Let $\pi:Y\rightarrow X$ as before, $G=\mathrm{Gal}(Y/X)$.\\
{\rm (i)} If $H\le G$ is not $\ell$-{\it hypoelementary},
\[ \sum _{U\le H}|U|\cdot \mu(U,H)\cdot m_{\ell}(\mathrm{Br}(Y/U)(\ell))=0. \]\\
{\rm (ii)} If $H\le G$ is not {\it hypoelementary},
\[ \sum _{U\le H}|U|\cdot \mu(U,H)\cdot m(\mathrm{Br}(Y/U))=0. \]
\end{cor}

\subsection{Example 1}

For a prime $\ell$ and an abelian group $A$, its corank is defined as $\mathrm{rank}_{\mathbb{Z}_{\ell}}(T_{\ell}(A))$, where $T_{\ell}(A)=\underset{\underset{n}{\longleftarrow}}{\lim}\ Ker(\ell^n:A\rightarrow A)$. Here we denote this by $\mathrm{rk}_{\ell}(A)$:
\[ \mathrm{rk}_{\ell}(A):=\mathrm{rank}_{\mathbb{Z}_{\ell}}(T_{\ell}(A)) \]

$\mathrm{Br}(X)(\ell)$ is known to be of finite corank, for example in the following cases (\cite{Grothendieck}):\\
- {\rm (C1)} $k$: a separably closed or finite field, $X$: of finite type /$k$, and proper or smooth $/k$, or $\text{char}(k)=0$ or $\dim X\le 2$.\\
- {\rm (C2)} $X$: of finite type $/Spec(\mathbb{Z})$, and smooth $/Spec(\mathbb{Z})$ or proper over $\exists$open $\subset Spec(\mathbb{Z})$.

Remark that if $Y/X$ is a finite \'{e}tale covering and if $X$ satisfies {\rm (C1)} or {\rm (C2)}, then so does $Y$.

\begin{exmp}
Assume $X$ satisfies {\rm (C1)} or {\rm (C2)}. If a subgroup $H\le G$ is not $\ell$-hypoelementary, we have an equation
\[ \sum_{U\le H}|U|\mu(U,H)\cdot \mathrm{rk}_{\ell}(\mathrm{Br}(Y/H)(\ell))=0. \]
\end{exmp}

\subsection{Example 2}
By Gabber's lemma (Lemma 4 in \cite{Gabber}), for any finite \'{e}tale covering $\pi:Y\rightarrow X$, we have
\[ Br^{\prime}(X)/\mathrm{Br}(X)\overset{\pi^{\ast}}{\hookrightarrow} Br^{\prime}(Y)/\mathrm{Br}(Y). \]

In particular, if $\mathrm{Br}(Y)\subset \mathrm{Br}(Y)^{\prime}$ is of finite index, 
then so is $\mathrm{Br}(X)\subset \mathrm{Br}(X)^{\prime}$.

\begin{exmp}
Assume $Y$ satisfies $[Br^{\prime}(Y):\mathrm{Br}(Y)]<\infty$. Then for any non-hypoelementary subgroup $H\le G$, we have an equation
\[ \sum_{U\le H}|U|\mu(U,H)\cdot [Br^{\prime}(Y/U):\mathrm{Br}(Y/U)]=0. \]
\end{exmp}

\section{Appendix 2}
We showed Theorem \ref{TheoremToShow+} by using \v{C}ech cohomology. This way of proof required Assumption \ref{Assumption}.

In this section, to get rid of Assumption \ref{Assumption}, 
we consider a more general proof.

\begin{rem}
By definition, $N_{\pi}:H^2_{\mathrm{et}}(Y,\mathbb{G}_{m,Y})\rightarrow H^2_{\mathrm{et}}(X,\mathbb{G}_{m,X})$ is the composition
\[ H^2_{\mathrm{et}}(Y,\mathbb{G}_{m,Y})
\overset{\mathfrak{c}^{-1}}{\underset{\cong}{\longrightarrow}}
H^2_{\mathrm{et}}(X,\pi_{\ast}\mathbb{G}_{m,Y})\overset{H^2_{\mathrm{et}}(X,N_{\pi})}{\longrightarrow}H^2_{\mathrm{et}}(X,\mathbb{G}_{m,X}), \]
where
\[ N_{\pi}:\pi_{\ast}\mathbb{G}_{m,Y}\rightarrow\mathbb{G}_{m,X} \]
is the norm homomorphism in $\mathbf{S}(X_{\mathrm{et}})$.

Thus the diagram in Theorem \ref{TheoremToShow+} is nothing other than the following $:$
\[
\xy
(-20,8)*+{\mathrm{Br}(Y)}="0";
(20,8)*+{\mathrm{Br}(X)}="2";
(-28,-4)*+{H^2_{\mathrm{et}}(Y,\mathbb{G}_{m,Y})}="4";
(0,-16)*+{H^2_{\mathrm{et}}(X,\pi_{\ast}\mathbb{G}_{m,Y})}="6";
(28,-4)*+{H^2_{\mathrm{et}}(X,\mathbb{G}_{m,X})}="8";
(0,8)*+{}="10";
{\ar^{N_{\pi}} "0";"2"};
{\ar_{\chi_Y} "0";"4"};
{\ar^{\chi_X} "2";"8"};
{\ar_{\mathfrak{c}^{-1}}^{\cong} "4";"6"};
{\ar_{H^2_{\mathrm{et}}(N_{\pi})} "6";"8"};
{\ar@{}|\circlearrowright "6";"10"};
\endxy
\]
Remark also that we may assume $X$ is connected.
\end{rem}

\begin{rem}\label{RemarkCompat1}
For any finite \'{e}tale covering $\pi:Y\rightarrow X$,
\[ \pi_{\ast}:\mathbf{S}(Y_{\mathrm{et}})\rightarrow\mathbf{S}(X_{\mathrm{et}}) \]
is exact. Here, $\mathbf{S}(X_{\mathrm{et}})$ denotes the category of abelian sheaves on $X_{\mathrm{et}}$. 
Thus we have natural homomorphisms
\[ \pi_{\ast}:H^q_{\mathrm{et}}(Y,\mathbb{G}_{m,Y})\rightarrow H^q_{\mathrm{et}}(X,\pi_{\ast}\mathbb{G}_{m,Y})\qquad(\forall q\ge0).  \]

It can be easily seen that this gives the inverse of
\[ \mathfrak{c}:H^2_{\mathrm{et}}(X,\pi_{\ast}\mathbb{G}_{m,Y})\rightarrow H^2_{\mathrm{et}}(Y,\mathbb{G}_{m,Y}). \]
\end{rem}

\begin{prop}
Let $\pi:Y\rightarrow X$ be a finite \'{e}tale covering of a connected scheme $X$. 
For any $\mathbb{G}_{m,Y}$-gerbe $F$ on $Y_{\mathrm{et}}$, if we define a fibered category $\pi_{\ast}F$ over $X_{\mathrm{et}}$ by
\[ (\pi_{\ast}F)(U)=F(U\times_XY)\qquad (\forall U\in X_{\mathrm{et}}), \]
in a natural way, then $\pi_{\ast}F$ becomes a $\pi_{\ast}\mathbb{G}_{m,Y}$-gerbe on $X_{\mathrm{et}}$.
This defines a group homomorphism
\[ \pi_{\ast}:H^2_g(Y,\mathbb{G}_{m,Y})\rightarrow H^2_g(X,\pi_{\ast}\mathbb{G}_{m,Y}), \]
where $H^2_g$ denotes the non-abelian cohomology of Giraud.

\end{prop}
\begin{pf}
Since $F$ is a stack fibered in groupoid, it can be easily seen that so is $\pi_{\ast}F$. 
Thus, to show $\pi_{\ast}F$ is a gerbe, it suffices to show the following:\\
\quad {\rm (a)} $\pi_{\ast}F$ is locally connected\\
\quad {\rm (b)} $\pi_{\ast}F$ is locally non-empty



{\rm (a)} For any $U\in X_{\mathrm{et}}$ and any $a_1,a_2\in\pi_{\ast}F(U)=F(V)$ $(V:=U\times_XY)$, there exists a surjective \'{e}tale morphism $V^{\prime}\overset{v}{\rightarrow}V$ of finite type 
such that $v^{\ast}a_1\cong v^{\ast}a_2$ in $F(V^{\prime})$.

By Remark \ref{ec}, there exists a surjective \'{e}tale morphism $U^{\prime}\overset{u}{\rightarrow}U$ of finite type 
such that $U^{\prime}\times_XY=U^{\prime}\times_UV\overset{\mathrm{pr}_Y}{\rightarrow}V$ factors through $v$:

\begin{equation}
\xy
(0,4)*+{V^{\prime}}="0";
(0,-10)*+{}="1";
(-24,-6)*+{U^{\prime}\times_UV}="2";
(24,-6)*+{V}="4";
{\ar^{\exists w} "2";"0"};
{\ar^{v} "0";"4"};
{\ar@/_0.80pc/_{\mathrm{pr}_V} "2";"4"};
{\ar@{}|\circlearrowright "0";"1"};
\endxy
\label{FactorDiag}
\end{equation}

Thus we have $w^{\ast}v^{\ast}a_1\cong w^{\ast}v^{\ast}a_2$ in $F(U^{\prime}\times_UV)$, namely, $u^{\ast}a_1\cong u^{\ast}a_2$ in $\pi_{\ast}F(U^{\prime})$.

{\rm (b)}
For any $U\in X_{\mathrm{et}}$, let $V^{\prime}\overset{v}{\rightarrow}V=U\times_XY$ be a surjective \'{e}tale morphism of finite type, such that $F(V^{\prime})\ne\emptyset$.
Take $U^{\prime}\overset{u}{\rightarrow}U$ satisfying $(\ref{FactorDiag})$ as above.

If we put $W_1:=w(U^{\prime}\times_UV)$ and $W_2:=V^{\prime}\setminus W_1$, then each $W_i$ is an open subscheme of $V^{\prime}$ $(i=1,2)$, satisfying
\[ V^{\prime}=W_1\coprod W_2. \]
Thus we have $F(V^{\prime})\simeq F(W_1)\times F(W_2)$. In particular, $F(W_1)\ne\emptyset$. Since $w:U^{\prime}\times_UV\rightarrow W_1$ is surjective \'{e}tale, $\pi_{\ast}F(U^{\prime})=F(U^{\prime}\times_UV)\ne\emptyset$ follows from $F(W_1)\ne\emptyset$.

Thus $\pi_{\ast}F$ is a gerbe, which is obviously bound by $\pi_{\ast}\mathbb{G}_{m,Y}$.
\end{pf}

\begin{rem}$(${\rm \cite{Milne}}$)$
Let $X$ be a scheme. For any Azumaya algebra $\mathcal{A}$ on $X$,
let $F_{\mathcal{A}}$ denote a
fibered category over $X_{\mathrm{et}}$,
whose fiber $F_{\mathcal{A}}(U)$ over $U\in X_{\mathrm{et}}$ is defined as follows $:$

\noindent- An object is a pair $(\mathcal{E},\alpha)$, where $\mathcal{E}\in\mathrm{q\text{-}Coh}(U)$ is locally free of finite rank, $\alpha:\mathcal{E}nd_{\mathcal{O}_U}(\mathcal{E})\overset{\cong}{\longrightarrow}\mathcal{A}\mid_U$ is an isomorphism of $\mathcal{O}_U$-algebras.

\noindent- A morphism $(\mathcal{E},\alpha)\rightarrow(\mathcal{E}^{\prime},\alpha^{\prime})$ is an isomorphism $\mathcal{E}\overset{\cong}{\rightarrow}\mathcal{E}^{\prime}$ compatible with $\alpha$ and $\alpha^{\prime}$.

Then $F_{\mathcal{A}}$ becomes a gerbe, bound by $\mathbb{G}_{m,X}$. $($Indeed, multiplication by elements of $\Gamma(U,\mathcal{O}_U^{\times})$ gives an isomorphism $
\Gamma(U,\mathcal{O}_U^{\times})\overset{\cong}{\rightarrow}\mathrm{Aut}_{F_{\mathcal{A}}}(\mathcal{E},\alpha)$.$)$

This defines the natural monomorphism
\[ \chi_X:\mathrm{Br}(X)\hookrightarrow H^2_g(X,\mathbb{G}_{m,X}). \]
\end{rem}

\begin{lem}\label{LemmaCompat4}
Let $\pi:Y\rightarrow X$ be a finite \'{e}tale covering, and let $\mathcal{A}$ be an Azumaya algebra on $Y$.

{\rm (i)}  For any $U\in X_{\mathrm{et}}$, let
\[
\xy
(-8,8)*+{V}="0";
(-8,-8)*+{U}="2";
(8,8)*+{Y}="4";
(8,-8)*+{X}="6";
(0,0)*+{\square}="10";
{\ar_{\varpi} "0";"2"};
{\ar^{\pi} "4";"6"};
{\ar^{g} "0";"4"};
{\ar_{f} "2";"6"};
\endxy
\]
be a pull-back diagram. We define a functor
\[ \mathcal{N}_{\varpi}:F_{\mathcal{A}}(V)\rightarrow F_{\mathcal{N}_{\pi}\mathcal{A}}(U) \]
by $\mathcal{N}_{\varpi}(\mathcal{E},\alpha)=(\mathcal{N}_{\varpi}(\mathcal{E}),\beta)$,
where $\beta$ is the composition
\[ \mathcal{E}nd_{\mathcal{O}_U}(\mathcal{N}_{\varpi}(\mathcal{E}))\overset{\cong}{\rightarrow}%
\mathcal{N}_{\varpi}(\mathcal{E}nd_{\mathcal{O}_V}(\mathcal{E}))\overset{\mathcal{N}_{\varpi}(\alpha)}{\longrightarrow}%
\mathcal{N}_{\varpi}(g^{\ast}\mathcal{A})\overset{\cong}{\rightarrow}%
f^{\ast}\mathcal{N}_{\pi}\mathcal{A}. \]

Then for any morphism $u:U^{\prime}\rightarrow U$ in $X_{\mathrm{et}}$ if we take the pull-back
\[
\xy
(-8,8)*+{V^{\prime}}="0";
(-8,-8)*+{U^{\prime}}="2";
(8,8)*+{V}="4";
(8,-8)*+{U}="6";
(11,-9)*+{,}="7";
(0,0)*+{\square}="10";
{\ar_{\varpi^{\prime}} "0";"2"};
{\ar^{\varpi} "4";"6"};
{\ar^{v} "0";"4"};
{\ar_{u} "2";"6"};
\endxy
\]
then we have a natural isomorphism
\[ u^{\ast}\mathcal{N}_{\varpi}\cong\mathcal{N}_{\varpi^{\prime}}v^{\ast}:F_{\mathcal{A}}(V)\rightarrow F_{\mathcal{N}_{\pi}\mathcal{A}}(U^{\prime}). \]

{\rm (ii)} $\mathcal{N}_{\varpi}$ makes the following diagram commutative $:$
\[
\xy
(-28,-8)*+{\mathrm{Aut}_{F_{\mathcal{A}}(V)}(\mathcal{E},\alpha)}="0";
(28,-8)*+{\mathrm{Aut}_{F_{\mathcal{N}_{\pi}\mathcal{A}}(U)}(\mathcal{N}_{\varpi}(\mathcal{E},\alpha))}="2";
(-28,8)*+{\Gamma(V,\mathcal{O}_V^{\times})}="4";
(28,8)*+{\Gamma(U,\mathcal{O}_U^{\times})}="6";
{\ar_{\mathcal{N}_{\varpi}} "0";"2"};
{\ar^{N_{\varpi}(U)} "4";"6"};
{\ar^{\cong} "4";"0"};
{\ar^{\cong} "6";"2"};
{\ar@{}|\circlearrowright "0";"6"};
\endxy
\]
\end{lem}
\begin{pf}
{\rm (i)} This is induced from the natural monoidal isomorphism
\[ \theta^u_{\varpi}:u^{\ast}\mathcal{N}_{\varpi}\overset{\cong}{\rightarrow}\mathcal{N}_{\varpi^{\prime}}v^{\ast}:\mathrm{q\text{-}Coh}(V)\rightarrow\mathrm{q\text{-}Coh}(U^{\prime}). \]

{\rm (ii)} This follows from the commutativity of
\[
\xy
(-16,-8)*+{\mathrm{Aut}_{\mathcal{O}_V}(\mathcal{E})}="0";
(16,-8)*+{\mathrm{Aut}_{\mathcal{O}_U}(\mathcal{N}_{\varpi}(\mathcal{E}))}="2";
(-16,8)*+{\Gamma(V,\mathcal{O}_V^{\times})}="4";
(16,8)*+{\Gamma(U,\mathcal{O}_U^{\times})}="6";
{\ar_{\mathcal{N}_{\varpi}} "0";"2"};
{\ar^{N_{\varpi}(U)} "4";"6"};
{\ar "4";"0"};
{\ar "6";"2"};
{\ar@{}|\circlearrowright "0";"6"};
\endxy
\text{and}
\xy
(-20,-8)*+{\mathrm{Aut}_{\mathcal{O}_V}(\mathcal{E})}="0";
(20,-8)*+{\mathrm{Aut}_{\mathcal{O}_U}(\mathcal{N}_{\varpi}(\mathcal{E}))}="2";
(-20,8)*+{\mathrm{Aut}_{F_{\mathcal{A}}(V)}(\mathcal{E},\alpha)}="4";
(20,8)*+{\mathrm{Aut}_{F_{\mathcal{N}_{\pi}\mathcal{A}}(U)}(\mathcal{N}_{\varpi}(\mathcal{E},\alpha))}="6";
(34,-9)*+{.}="7";
{\ar_{\mathcal{N}_{\varpi}} "0";"2"};
{\ar^<<<<{\mathcal{N}_{\varpi}} "4";"6"};
{\ar@{^(->} "4";"0"};
{\ar@{^(->} "6";"2"};
{\ar@{}|\circlearrowright "0";"6"};
\endxy
\]
\end{pf}

\begin{rem}$(${\rm \cite{Giraud}} Proposition 3.1.5$)$
Let $X$ be a scheme. For any morphism $u:\mathcal{F}\rightarrow\mathcal{G}$ in $\mathbf{S}(X_{\mathrm{et}})$, we have a group homomorphism
\[ H^2_g(u):H^2_g(X,\mathcal{F})\rightarrow H^2_g(X,\mathcal{G}). \] 
If $F$ is an $\mathcal{F}$-gerbe and $G$ is a $\mathcal{G}$-gerbe, then $H^2_g(u)(F)=G$ in $H^2_g(X,\mathcal{G})$ if and only if there exists a morphism of gerbes $F\rightarrow G$ bound by $u$.
\end{rem}

By the above arguments, Theorem \ref{TheoremToShow+} is reduced to the following Proposition:
\begin{prop}\label{MainLemCompat}
Let $\pi:Y\rightarrow X$ be a finite \'{e}tale covering of a connected scheme $X$. The following diagram is commutative:
\[
\xy
(-20,8)*+{\mathrm{Br}(Y)}="0";
(20,8)*+{\mathrm{Br}(X)}="2";
(-28,-4)*+{H^2_g(Y,\mathbb{G}_{m,Y})}="4";
(0,-16)*+{H^2_g(X,\pi_{\ast}\mathbb{G}_{m,Y})}="6";
(28,-4)*+{H^2_g(X,\mathbb{G}_{m,X})}="8";
(0,8)*+{}="10";
{\ar^{N_{\pi}} "0";"2"};
{\ar_{\chi_Y} "0";"4"};
{\ar^{\chi_X} "2";"8"};
{\ar_{\pi_{\ast}} "4";"6"};
{\ar_{H^2_g(N_{\pi})} "6";"8"};
{\ar@{}|\circlearrowright "6";"10"};
\endxy
\]
\end{prop}
\begin{pf}
By {\rm (i)} in Lemma \ref{LemmaCompat4},
if we attach a functor
\[ \mathcal{N}_{\pi}(U):=\mathcal{N}_{\varpi}:\pi_{\ast}F_{\mathcal{A}}(U)\rightarrow F_{\mathcal{N}_{\pi}\mathcal{A}}(U) \]
to each $U\in X_{\mathrm{et}}$,
then
\[ \mathcal{N}_{\pi}(U)\qquad(U\in X_{\mathrm{et}}) \]
forms a morphism of fibered categories
\[ \mathcal{N}_{\pi}:\pi_{\ast}F_{\mathcal{A}}\rightarrow F_{\mathcal{N}_{\pi}\mathcal{A}}. \]

By {\rm (ii)} in Lemma \ref{LemmaCompat4}, this is bound by $N_{\pi}:\pi_{\ast}\mathbb{G}_{m,Y}\rightarrow\mathbb{G}_{m,X}$.

Thus we have $H^2_g(N_{\pi})(\pi_{\ast}F_{\mathcal{A}})=F_{\mathcal{N}_{\pi}\mathcal{A}}$.
\end{pf}




\begin{thebibliography}{00}




\bibitem{Bley-Boltje}
W. Bley, R. Boltje, Cohomological Mackey functors in number theory, J. Number Theory \textbf{105} (2004), 1-37.

\bibitem{Boltje}
R. Boltje, Class group relations from Burnside ring idempotents, J. Number Theory \textbf{66} (1997), 291-305.

\bibitem{Bouc}
S. Bouc, Green functors and $G$-sets, Lecture Notes in Mathematics \textbf{1671} Springer-Verlag, Berlin, 1997. viii+342 pp.

\bibitem{Ferrand}
D. Ferrand, Un foncteur norme, Bull. Soc. math. France \textbf{126} (1998), 1-49.

\bibitem{Ford}
T. J. Ford, Hecke actions on Brauer groups, J. Pure Appl. Algebra \textbf{33} (1984), 11-17.

\bibitem{Gabber}
O. Gabber, Some theorems on Azumaya algebras, Lecture Notes in Math \textbf{844}  Springer-Verlag (1980), 129-210.

\bibitem{Giraud}
J. Giraud, Cohomologie non ab\'{e}lienne, Die Grundlehren der Mathematischen Wissenschaften in Einzeldarstellungen, Band \textbf{179}, Springer-Verlag, (1971).

\bibitem{Grothendieck}
A. Grothendieck, Le groupe de Brauer I,II,III, In Dix Expos\'{e}s sur la Cohomologie des Sch\'{e}mas, North-Holland, Amsterdam, (1968), 46-188.

\bibitem{Knus-Ojanguren}
M.-A. Knus, M. Ojanguren, A Norm for modules and algebras, Math. Z. \textbf{142} (1975), 33-45.

\bibitem{Milne}
J. S. Milne, \'{E}tale cohomology, Princeton Mathematical Series \textbf{33} Princeton University Press, Princeton, N.J., 1980. xiii+323 pp.

\bibitem{Murre}
J. P. Murre, Lectures on an introduction to Grothendieck's theory of the fundamental group, Notes by S. Anantharaman, Tata Institute of Fundamental Research Lectures on Mathmatics, No 40. Tata Institute of Fundamental Research, Bombay, 1967. iv+176+iv pp.


\end{thebibliography}
\end{document}